\documentclass[11pt,english]{article}
\usepackage[T1]{fontenc}
\usepackage[latin9]{inputenc}
\usepackage{geometry}
\usepackage{color}
\usepackage{url}
\usepackage{amsmath}
\usepackage{amssymb}

\makeatletter

\newcommand*\LyXZeroWidthSpace{\hspace{0pt}}

\usepackage{babel}

\makeatother

\usepackage{babel}
\begin{document}

\title{On the Bhattacharya-Mesner rank of third order hypermatrices}

\author{Edinah K. Gnang \thanks{Department of Applied Mathematics and Statistics, Johns Hopkins University,
email: egnang1@jhu.edu}, Yuval Filmus \thanks{Computer Science Department, Technion - Israel Institute of Technology}}
\maketitle
\begin{abstract}
We introduce the Bhattacharya-Mesner rank of third order hypermatrices
as a relaxation to the tensor rank and devise from it some bounds
for the tensor rank. We use the Bhattacharya-Mesner rank to extend
to third order hypermatrices the connection relating the rank to a
notion of linear dependence. We also derive explicit necessary and
sufficient conditions for the existence of third order hypermatrix
inverse pair. Finally we use inverse pair to extend to third order
hypermatrices the formulation and proof of the matrix rank\textendash nullity
theorem.
\end{abstract}

\section{Introduction}

\emph{Hypermatrices} are multidimensional analog of matrices. A hypermatrix
is therefore a finite multiset whose elements ( called entries ) are
indexed by distinct elements of some fixed Cartesian product of the
form
\[
\left\{ 0,\cdots,n_{0}-1\right\} \times\left\{ 0,\cdots,n_{1}-1\right\} \times\cdots\times\left\{ 0,\cdots,n_{m-1}-1\right\} .
\]
Such a hypermatrix is of order $m$ and of size $n_{0}\times n_{1}\times\cdots\times n_{m-1}$.
A hypermatrix is cubic, of side length $n$ if $n_{0}=n_{1}=\cdots=n_{m-1}=n$.
In particular, matrices correspond to second order hypermatrices.
Hypermatrix algebras arise from generalizations of classical matrix
notions and algorithms \cite{MB94,GKZ,RK,GER,MB90,GNANG2017238}.
The distinction between hypermatrices and tensors closely mirrors
the distinction between matrices and abstract linear transformations.
Recall that an abstract linear transformation defined over a finite
dimensional vector space is identified with a matrix orbit, not just
a single matrix. The orbit correspond to the different choices of
bases for the domain and range of the linear transformation. For instance,
let $\mathbf{M}\in\mathbb{K}^{m\times n}$ denote the matrix associated
with some abstract linear transformation relative to the standard
basis for $\mathbb{K}^{m\times1}$ and $\mathbb{K}^{1\times n}$.
Then the tensorial orbit associated with the corresponding abstract
linear transformation is defined to be the set of matrices given by
\[
\left\{ \mathbf{A}\cdot\mathbf{M}\cdot\mathbf{B}\,:\,\begin{array}{c}
\mathbf{A}\in\text{GL}_{m}\left(\mathbb{K}\right)\\
\mathbf{B}\in\text{GL}_{n}\left(\mathbb{K}\right)
\end{array}\right\} .
\]
A matrix property common to all the members of the same tensorial
orbit is called a \emph{tensorial invariant}. Classical matrix attributes
well known to be tensorial invariants when $\mathbb{K}=\mathbb{C}$
are: 
\begin{itemize}
\item The \emph{rank} defined to be 
\[
\min_{\begin{array}{c}
\mathbf{A}\in\text{GL}_{m}\left(\mathbb{K}\right)\\
\mathbf{B}\in\text{GL}_{n}\left(\mathbb{K}\right)
\end{array}}\left\Vert \mathbf{A}\cdot\mathbf{M}\cdot\mathbf{B}\right\Vert _{\ell_{0}}
\]
\item The \emph{nullity} defined to be the dimension of the nullspaces of
$\mathbf{A}$ and well known to be given by 
\[
\min\left(m,n\right)-\min_{\begin{array}{c}
\mathbf{A}\in\text{GL}_{m}\left(\mathbb{K}\right)\\
\mathbf{B}\in\text{GL}_{n}\left(\mathbb{K}\right)
\end{array}}\left\Vert \mathbf{A}\cdot\mathbf{M}\cdot\mathbf{B}\right\Vert _{\ell_{0}}
\]
\item Singular values defined to be the modulus of numbers in the multiset
formed by taking diagonal entries of any diagonal matrix element of
the tensorial sub-orbit restricted to left and right action by the
unitary subgroup of the general linear groups described as follows:
\[
\left\{ \mathbf{A}\cdot\mathbf{M}\cdot\mathbf{B}\,:\,\begin{array}{c}
\mathbf{A}\in\text{U}_{m}\left(\mathbb{K}\right)\\
\mathbf{B}\in\text{U}_{n}\left(\mathbb{K}\right)
\end{array}\right\} .
\]
\item The eigenvalues of square diagonalizable matrices defined to be the
multiset formed by taking the diagonal entries of any diagonal matrix
elements of the tensorial sub-orbit restricted to the action by conjugation
of the general linear group described as follows:
\[
\left\{ \mathbf{A}\cdot\mathbf{T}\cdot\mathbf{B}\,:\,\mathbf{A},\,\mathbf{B}\in\text{GL}_{n}\left(\mathbb{K}\right)\;\text{ and }\;\mathbf{A}\cdot\mathbf{B}=\mathbf{I}_{n}\right\} .
\]
\end{itemize}
Classically, hypermatrices such as third order hypermatrices from
$\mathbb{K}^{m\times n\times p}$ also arise from tensorial orbits
induced by the action of various appropriate subgroups of the general
linear group on canonical embeddings of the vector spaces $\mathbb{K}^{m\times1\times1}$,
$\mathbb{K}^{1\times n\times1}$ and $\mathbb{K}^{1\times1\times p}$
respectively. Incidentally, the classical tensor rank and singular
values are defined by analogy to their matrix counterparts. Unfortunately,
the tensorial invariant approach to defining the eigenvalues does
not extend to odd order hypermatrices because we cannot express the
action by conjugation using a triplet of matrix elements from the
general linear group. Fortunately the Bhattacharya-Mesner algebra
suggest new tensorial hypermatrix orbit based upon higher order analog
of the general linear group. Hypermatrix analog of general linear
group allows us to extend to hypermatrices the nottion of a spectral
decomposition \cite{GER,GNANG2017238}. We argue in the present note
that the Bhattacharya-Mesner (BM for short) rank also preserves the
link relating the rank to the nullity of third order hypermatrices
when these notions are appropriately defined. Our discussion therefore
focuses on the Bhattacharya-Mesner algebra first proposed in \cite{MB90,MB94},
and the general Bhattacharya-Mesner algebra first proposed in \cite{GER}.
The general Bhattacharya-Mesner product encompasses as special cases
many other hypermatrix products and decompositions discussed in the
literature, including the usual matrix product, the Segre outer product,
the contraction product, the higher order singular value decomposition,
the multilinear matrix multiplication \cite{Lim2013}, and the slice
rank factorization \cite{2016arXiv160506702B}. Note that the BM algebra
of third order hypermatrices in particular is also motivated by the
study of linear systems of equations defined over a skew field ( such
as the skew field of quaternions ) whose unknowns are multiplied both
on the left and on the right by known coefficients from the skew field.

This article is accompanied by an extensive and actively maintained
Sage \cite{sage} symbolic hypermatrix algebra package which implements
various features of the general Bhattacharya-Mesner algebra. The package
is made available at the link \url{https://github.com/gnang/HypermatrixAlgebraPackage}.\\
\\
\emph{Acknowledgement}: We would like to thank Andrei Gabrielov for
providing guidance while preparing this manuscript. We are grateful
to Dan Naiman for pointing out an upper bound on the diagonal dependence.
The first author was supported by the National Science Foundation
grant DMS\textendash 1161629\LyXZeroWidthSpace , and is grateful for
the hospitality of the Institute for Advanced Study and the Department
of Mathematics at Purdue University.

\section{Overview of the Bhattacharya-Mesner algebra of third order hypermatrices}

The Bhattacharya-Mesner product, first proposed in \cite{MB90,MB94},
generalizes the classical matrix product. The Bhattacharya-Mesner
product is best introduced by emphasizing its similarities with the
matrix product. Recall that for conformable matrices with entries
from some skew field $\mathbb{K}$
\[
\mathbf{A}^{(0)}\in\mathbb{K}^{n_{0}\times\textcolor{red}{\ell}}\,\mbox{ and }\mathbf{A}^{(1)}\in\mathbb{K}^{\textcolor{red}{\ell}\times n_{1}},
\]
we denote their product by $\mbox{Prod}\left(\mathbf{A}^{(0)},\mathbf{A}^{(1)}\right)\in\mathbb{K}^{n_{0}\times n_{1}}$
with entries given by 
\[
\mbox{Prod}\left(\mathbf{A}^{(0)},\mathbf{A}^{(1)}\right)\left[i_{0},i_{1}\right]=\sum_{0\le\textcolor{red}{j}<\ell}\mathbf{A}^{(0)}\left[i_{0},\textcolor{red}{j}\right]\mathbf{A}^{(1)}\left[\textcolor{red}{j},i_{1}\right].
\]
Whereas the matrix product is a binary operation, the BM product of
third order hypermatrices is a ternary operation. The product of third
order conformable hypermatrices
\[
\mathbf{A}^{(0)}\in\mathbb{K}^{n_{0}\times\textcolor{red}{\ell}\times n_{2}},\,\mathbf{A}^{(1)}\in\mathbb{K}^{n_{0}\times n_{1}\times\textcolor{red}{\ell}}\,\mbox{ and }\mathbf{A}^{(2)}\in\mathbb{K}^{\textcolor{red}{\ell}\times n_{1}\times n_{2}},
\]
denoted Prod$\left(\mathbf{A}^{(0)},\mathbf{A}^{(1)},\mathbf{A}^{(2)}\right)\in\mathbb{K}^{n_{0}\times n_{1}\times n_{2}}$
is specified entry-wise by
\[
\mbox{Prod}\left(\mathbf{A}^{(0)},\mathbf{A}^{(1)},\mathbf{A}^{(2)}\right)\left[i_{0},i_{1},i_{2}\right]=\sum_{0\le\textcolor{red}{j}<\ell}\mathbf{A}^{(0)}\left[i_{0},\textcolor{red}{j},i_{2}\right]\mathbf{A}^{(1)}\left[i_{0},i_{1},\textcolor{red}{j}\right]\mathbf{A}^{(2)}\left[\textcolor{red}{j},i_{1},i_{2}\right].
\]
In hindsight, the BM product of third order hypermatrices occurs naturally
in descriptions of \emph{general linear systems} of equations over
skew fields. A general linear systems of equations over a skew field
$\mathbb{K}$, is one for which the unknowns are multiplied both on
the left and on the right by coefficients from $\mathbb{K}$. General
linear system of equations are therefore of the form 
\begin{equation}
\left\{ \sum_{0\le\textcolor{red}{j}<\ell}\mathbf{A}\left[0,\textcolor{red}{j},i_{2}\right]\mathbf{x}\left[0,0,\textcolor{red}{j}\right]\mathbf{B}\left[\textcolor{red}{j},0,i_{2}\right]\,=\,\mathbf{c}\left[0,0,i_{2}\right]\right\} _{0\le i_{2}<n_{2}}.\label{Set of General Constraints}
\end{equation}
The Bhattacharya-Mesner product succinctly expresses such systems
in terms of a left-coefficient hypermatrix $\mathbf{A}\in\mathbb{K}^{1\times\textcolor{red}{\ell}\times n_{2}}$
, an unknown vector $\mathbf{x}\in\mathbb{K}^{1\times1\times\textcolor{red}{\ell}}$,
a right-coefficient hypermatrix $\mathbf{B}\in\mathbb{K}^{\textcolor{red}{\ell}\times1\times n_{2}}$
and a right hand side vector $\mathbf{c}\in\mathbb{K}^{1\times1\times n_{2}}$
by the equation
\begin{equation}
\mbox{Prod}\left(\mathbf{A},\,\mathbf{x},\,\mathbf{B}\right)=\mathbf{c}.\label{General Linear Constraints}
\end{equation}
Note that the left hand side of a general linear system of equations
defined over skew fields as described in Eq. (\ref{Set of General Constraints})
can not generally be expressed as a single left action ( or right
action ) of some coefficient matrix on an vector of unknowns. Special
instances of general linear systems of equations can be solved using
quasi-determinants \cite{GR,GRW,Gelfand1997}, but in general solutions
to such systems are not expressible as non-commutative rational functions
of the entries of $\mathbf{A}$, $\mathbf{B}$ and $\mathbf{c}$.

We further recall a variant of the Bhattacharya-Mesner product called
the \emph{general Bhattacharya-Mesner product} of particular interest
to our discussion because it simplifies subsequent notations. The
general product of third order hypermatrices is defined for any conformable
triple 
\[
\mathbf{A}^{(0)}\in\mathbb{K}^{n_{0}\times\textcolor{red}{\ell}\times n_{2}},\,\mathbf{A}^{(1)}\in\mathbb{K}^{n_{0}\times n_{1}\times\textcolor{red}{\ell}}\,\mbox{ and }\mathbf{A}^{(2)}\in\mathbb{K}^{\textcolor{red}{\ell}\times n_{1}\times n_{2}},
\]
taken with an additional cubic \emph{background hypermatrix} $\mathbf{B}$
of side length $\textcolor{red}{\ell}$ ( the contracted dimension
). The general Bhattacharya-Mesner product of hypermatrices $\mathbf{A}^{(0)}$,
$\mathbf{A}^{(1)}$ and $\mathbf{A}^{(2)}$ taken with the background
hypermatrix $\mathbf{B}$, denoted $\mbox{Prod}_{\mathbf{B}}\left(\mathbf{A}^{(0)},\mathbf{A}^{(1)},\mathbf{A}^{(2)}\right)\in\mathbb{K}^{n_{0}\times n_{1}\times n_{2}}$
is specified entry-wise by
\[
\mbox{Prod}_{\mathbf{B}}\left(\mathbf{A}^{(0)},\mathbf{A}^{(1)},\mathbf{A}^{(2)}\right)\left[i_{0},i_{1},i_{2}\right]=
\]
\begin{equation}
\sum_{0\le\textcolor{red}{j_{0}},\textcolor{red}{j_{1}},\textcolor{red}{j_{2}}<\ell}\mathbf{A}^{(0)}\left[i_{0},\textcolor{red}{j_{0}},i_{2}\right]\mathbf{A}^{(1)}\left[i_{0},i_{1},\textcolor{red}{j_{1}}\right]\mathbf{A}^{(2)}\left[\textcolor{red}{j_{2}},i_{1},i_{2}\right]\mathbf{B}\left[\textcolor{red}{j_{0}},\textcolor{red}{j_{1}},\textcolor{red}{j_{2}}\right].\label{General BM product}
\end{equation}

The original Bhattacharya-Mesner product is recovered from the general
product by setting the cubic background hypermatrix $\mathbf{B}$
to be equal to the \emph{Kronecker delta} hypermatrix denoted $\boldsymbol{\Delta}$,
whose entries are specified by 
\[
\boldsymbol{\Delta}\left[i_{0},i_{1},i_{2}\right]=\begin{cases}
\begin{array}{cc}
1 & \mbox{ if }\:0\le i_{0}=i_{1}=i_{2}<n\\
0 & \mbox{otherwise}
\end{array}\end{cases}.
\]

Finally we recall for the reader's convenience the definition of the
hypermatrix transpose operations. Let $\mathbf{A}\in\mathbb{C}^{n_{0}\times n_{1}\times n_{2}}$,
the corresponding transpose, denoted $\mathbf{A}^{\top}\in\mathbb{C}^{n_{1}\times n_{2}\times n_{0}}$,
is specified entry-wise by
\[
\mathbf{A}^{\top}\left[i_{0},i_{1},i_{2}\right]=\mathbf{A}\left[i_{2},i_{0},i_{1}\right].
\]
The transpose operation therefore performs a cyclic permutation to
the indices of each entry. For notational convenience we adopt the
convention
\[
\mathbf{A}^{\top^{2}}:=\left(\mathbf{A}^{\top}\right)^{\top},\ \mathbf{A}^{\top^{3}}:=\left(\mathbf{A}^{\top^{2}}\right)^{\top}=\mathbf{A}.
\]
By this convention, 
\[
\mathbf{A}^{\top^{i}}=\mathbf{A}^{\top^{j}}\:\mbox{ if }\:i\equiv j\mod3.
\]

\section{The Bhattacharya-Mesner rank}

The Bhattacharya-Mesner\emph{ outer-product} corresponds to the special
product instance where the the contracted dimension equals $1$. Furthermore,
the general Bhattacharya-Mesner product provides a convenient alternative
way of expressing BM outer-products. Recall that for an arbitrary
conformable triple
\[
\mathbf{A}^{(0)}\in\mathbb{K}^{n_{0}\times\textcolor{red}{\ell}\times n_{2}},\,\mathbf{A}^{(1)}\in\mathbb{K}^{n_{0}\times n_{1}\times\textcolor{red}{\ell}}\,\mbox{ and }\mathbf{A}^{(2)}\in\mathbb{K}^{\textcolor{red}{\ell}\times n_{1}\times n_{2}},
\]
an outer product corresponds to a product of selected oriented matrix
slices
\begin{equation}
\mbox{Prod}\left(\mathbf{A}^{(0)}\left[:,t,:\right],\,\mathbf{A}^{(1)}\left[:,:,t\right],\,\mathbf{A}^{(2)}\left[t,:,:\right]\right).\label{Outerproduct}
\end{equation}
Recall that in the colon notation, $\mathbf{A}^{(0)}\left[:,t,:\right]$
refers to a hypermatrix \emph{slice} of size $n_{0}\times1\times n_{2}$,
where the second index is fixed to $t$ while all other indices are
allowed to vary within their ranges. The BM product and the usual
matrix product have in common the fact that every product is a sum
of outer products. Hypermatrix outer products are conveniently expressible
as general BM products. The corresponding background hypermatrices
are denoted $\left\{ \boldsymbol{\Delta}^{(t)}\right\} _{0\le t<n}$
and specified entry-wise by 
\[
\boldsymbol{\Delta}^{(t)}\left[i_{0},i_{1},i_{2}\right]=\begin{cases}
\begin{array}{cc}
1 & \mbox{ if }\:0\le t=i_{0}=i_{1}=i_{2}<n\\
0 & \mbox{otherwise}
\end{array}\end{cases}.
\]
The outer product in Eq. (\eqref{Outerproduct}) is thus equal to
$\mbox{Prod}_{\boldsymbol{\Delta}^{(t)}}\left(\mathbf{A}^{(0)},\mathbf{A}^{(1)},\mathbf{A}^{(2)}\right)$.
The Bhattacharya-Mesner outer product suggests a natural generalization
of the notion of rank which differs from the classical tensor rank.
To emphasize the similarities with the matrix rank, recall that a
matrix $\mathbf{A}\in\mathbb{K}^{n_{0}\times n_{1}}$ has rank $r$
( over $\mathbb{K}$ ) if there exists a conformable matrix pair
\[
\mathbf{X}^{(0)}\in\mathbb{K}^{n_{0}\times\textcolor{red}{\ell}}\,\mbox{ and }\mathbf{X}^{(1)}\in\mathbb{K}^{\textcolor{red}{\ell}\times n_{1}},
\]
 such that 
\[
\mathbf{A}=\sum_{0\le t<r}\mbox{Prod}_{\boldsymbol{\Delta}^{(t)}}\left(\mathbf{X}^{(0)},\mathbf{X}^{(1)}\right),
\]
and crucially, for all conformable matrix pair 
\[
\mathbf{Y}^{(0)}\in\mathbb{K}^{n_{0}\times\textcolor{red}{\ell}}\,\mbox{ and }\mathbf{Y}^{(1)}\in\mathbb{K}^{\textcolor{red}{\ell}\times n_{1}},
\]
we have 
\[
\mathbf{A}\ne\sum_{0\le t<r-1}\mbox{Prod}_{\boldsymbol{\Delta}^{(t)}}\left(\mathbf{Y}^{(1)},\mathbf{Y}^{(2)}\right),
\]
where 
\[
\boldsymbol{\Delta}^{(t)}\left[i_{0},i_{1}\right]=\begin{cases}
\begin{array}{cc}
1 & \mbox{ if }\:0\le t=i_{0}=i_{1}<n\\
0 & \mbox{otherwise}
\end{array}\end{cases}.
\]
In other words, the matrix rank is the minimum number of outer products
which add up to $\mathbf{A}$. This approach to defining the rank
extends verbatim to hypermatrices and is called the Bhattacharya-Mesner
rank of a hypermatrix. Throughout our discussion unless otherwise
specified the rank refers to the Bhattacharya-Mesner rank. A hypermatrix
$\mathbf{A}\in\mathbb{K}^{n_{0}\times n_{1}\times n_{2}}$ has rank
$r$ (over $\mathbb{K}$) if there exists a conformable triple 
\[
\mathbf{X}^{(0)}\in\mathbb{K}^{n_{0}\times\textcolor{red}{\ell}\times n_{2}},\,\mathbf{X}^{(1)}\in\mathbb{K}^{n_{0}\times n_{1}\times\textcolor{red}{\ell}}\,\mbox{ and }\mathbf{X}^{(2)}\in\mathbb{K}^{\textcolor{red}{\ell}\times n_{1}\times n_{2}},
\]
 such that 
\begin{equation}
\mathbf{A}=\sum_{0\le t<r}\mbox{Prod}_{\boldsymbol{\Delta}^{(t)}}\left(\mathbf{X}^{(0)},\mathbf{X}^{(1)},\mathbf{X}^{(2)}\right),\label{Decomposition}
\end{equation}
and crucially for all BM conformable triple 
\[
\mathbf{Y}^{(0)}\in\mathbb{K}^{n_{0}\times\textcolor{red}{\ell}\times n_{2}},\,\mathbf{Y}^{(1)}\in\mathbb{K}^{n_{0}\times n_{1}\times\textcolor{red}{\ell}}\,\mbox{ and }\mathbf{Y}^{(2)}\in\mathbb{K}^{\textcolor{red}{\ell}\times n_{1}\times n_{2}},
\]
 we have 
\[
\mathbf{A}\ne\sum_{0\le t<r-1}\mbox{Prod}_{\boldsymbol{\Delta}^{(t)}}\left(\mathbf{Y}^{(0)},\mathbf{Y}^{(1)},\mathbf{Y}^{(2)}\right),
\]
where 
\[
\boldsymbol{\Delta}^{(t)}\left[i_{0},i_{1},i_{2}\right]=\begin{cases}
\begin{array}{cc}
1 & \mbox{ if }\:0\le t=i_{0}=i_{1}=i_{2}<n\\
0 & \mbox{otherwise}
\end{array}\end{cases}.
\]
In other words, the rank of $\mathbf{A}$ is the minimum number of
outer products which add up to $\mathbf{A}$. Note that the usual
notions of tensor/hypermatrix rank discussed in the literature \cite{Lim2013},
including the canonical polyadic rank, and more recently the slice
rank \cite{2016arXiv160506702B} are all special instances of the
Bhattacharya-Mesner rank where additional constraints are imposed
on the hypermatrices in the conformable triple $\left(\mathbf{X}^{(0)},\mathbf{X}^{(1)},\mathbf{X}^{(2)}\right)$
in Eq. (\ref{Decomposition}), as illustrated by the following proposition.\\
\\
\textbf{Proposition 1} : 
\[
\text{BM-rank}\left(\mathbf{A}\right)\le\text{tensor rank}\left(\mathbf{A}\right),\quad\forall\:\mathbf{A}\in\mathbb{K}^{m\times n\times p}.
\]
\\
\emph{Proof} : Recall the definition of the Kronecker product of vectors
denoted $\left(\mathbf{x}\otimes\mathbf{y}\otimes\mathbf{z}\right)\in\mathbb{K}^{m\times n\times p}$
of vectors $\mathbf{x}\in\mathbb{K}^{m\times1\times1}$, $\mathbf{y}\in\mathbb{K}^{1\times n\times1}$,
$\mathbf{z}\in\mathbb{K}^{1\times1\times p}$ associated with the
canonical decomposition specified entry wise by 
\[
\left(\mathbf{x}\otimes\mathbf{y}\otimes\mathbf{z}\right)\left[i,j,k\right]=\mathbf{x}\left[i,0,0\right]\cdot\mathbf{y}\left[0,j,0\right]\cdot\mathbf{z}\left[0,0,k\right].
\]
We justify the upper bound claim by establishing that vector outer
product above is a constrained Bhattacharya-Mesner outer product.
Let $\mathbf{X}\in\mathbb{C}^{m\times1\times p}$, $\mathbf{Y}\in\mathbb{C}^{m\times n\times1}$
and $\mathbf{Z}\in\mathbb{C}^{1\times n\times p}$ such that 
\[
\forall\,0\le k<p,\quad\mathbf{X}\left[i,0,k\right]=\mathbf{x}\left[i,0,0\right]=\mathbf{X}\left[i,0,k\right],
\]
\[
\forall\,0\le i<m,\quad\mathbf{Y}\left[i,j,0\right]=\mathbf{y}\left[0,j,0\right]=\mathbf{Y}\left[i,j,0\right],
\]
\[
\forall\,0\le j<n,\quad\mathbf{Z}\left[0,j,k\right]=\mathbf{z}\left[0,0,k\right]=\mathbf{Z}\left[0,j,k\right].
\]
We see that 
\[
\mathbf{x}\otimes\mathbf{y}\otimes\mathbf{z}=\mbox{Prod}_{\boldsymbol{\Delta}^{(0)}}\left(\mathbf{X},\mathbf{Y},\mathbf{Z}\right)=\mbox{Prod}\left(\mathbf{X},\mathbf{Y},\mathbf{Z}\right).
\]
The gap between the tensor rank and the Bhattacharya-Mesner rank is
well illustrated by the fact that the hypermatrix 
\[
\sum_{0\le t<r}\boldsymbol{\Delta}^{(t)}
\]
has tensor rank $r$ while its Bhattacharya-Mesner rank is $1$ for
all $0<r\le n$. So that the gap between the canonical polyadic rank
and the Bhattacharya-Mesner rank can be as large as $n-1$ for some
$n\times n\times n$ hypermatrix. A similar argument also establishes
the slice rank discussed in \cite{2016arXiv160506702B} as an upper
bound to the Bhattacharya-Mesner rank.\\
\\
The following proposition emphasize the similarity of between the
Bhattacharya-Mesner rank and the matrix rank.\\
\\
\\
\textbf{Proposition 2} : 
\[
\text{BM-rank}\left(\mathbf{A}\right)\le\min\left(m,\,n,\,p\right),\quad\forall\:\mathbf{A}\in\mathbb{K}^{m\times n\times p}.
\]
The Bhattacharya-Mesner rank of $\mathbf{A}\in\mathbb{C}^{m\times n\times p}$
is upper bounded by min$\left(m,n,p\right)$.\\
\\
\\
\emph{Proof} : Let $\mathbf{A}\in\mathbb{K}^{m\times n\times p}$
where $p\le\min\left\{ m,n\right\} $. Let $\mathbf{J}_{0}\in\mathbb{K}^{m\times p\times p}$
and $\mathbf{J}_{1}\in\mathbb{C}^{p\times n\times p}$ with entries
specified by 
\[
\mathbf{J}_{0}\left[i,t,k\right]=\begin{cases}
\begin{array}{cc}
1 & \mbox{ if }0\le t=k<p\\
0 & \mbox{otherwise}
\end{array} & \forall\,\begin{cases}
\begin{array}{c}
0\le i<m\\
0\le t<p\\
0\le k<p
\end{array}\end{cases}\end{cases},
\]
\[
\mathbf{J}_{1}\left[t,j,k\right]=\begin{cases}
\begin{array}{cc}
1 & \mbox{ if }0\le t=k<p\\
0 & \mbox{otherwise}
\end{array} & \forall\,\begin{cases}
\begin{array}{c}
0\le t<p\\
0\le j<n\\
0\le k<p
\end{array}\end{cases}\end{cases},
\]
consequently we have 
\[
\mathbf{A}=\mbox{Prod}\left(\mathbf{J}_{0},\mathbf{A},\mathbf{J}_{1}\right)=\sum_{0\le t<p}\mbox{Prod}_{\boldsymbol{\Delta}^{(t)}}\left(\mathbf{J}_{0},\mathbf{A},\mathbf{J}_{1}\right)
\]
which establishes that the rank of $\mathbf{A}$ must be less or equal
to $p$. Recall the transpose identity 
\[
\mathbf{A}^{\top}=\mbox{Prod}\left(\mathbf{A}^{\top},\mathbf{J}_{1}^{\top},\mathbf{J}_{0}^{\top}\right)=\sum_{0\le t<p}\mbox{Prod}_{\boldsymbol{\Delta}^{(t)}}\left(\mathbf{A}^{\top},\mathbf{J}_{1}^{\top},\mathbf{J}_{0}^{\top}\right)
\]
and 
\[
\mathbf{A}^{\top^{2}}=\mbox{Prod}\left(\mathbf{J}_{1}^{\top^{2}},\mathbf{J}_{0}^{\top^{2}},\mathbf{A}^{\top^{2}}\right)=\sum_{0\le t<p}\mbox{Prod}_{\boldsymbol{\Delta}^{(t)}}\left(\mathbf{J}_{1}^{\top^{2}},\mathbf{J}_{0}^{\top^{2}},\mathbf{A}^{\top^{2}}\right)
\]
from which the desired result follows

\subsection{Rank inequalities}

We introduce here the notion of matrix diagonal independence which
is closely related to the rank of third order hypermatrices. Let 
\[
\mathbf{A}\in\mathbb{K}^{m\times\ell\times p},\ \mathbf{B}\in\mathbb{K}^{m\times n\times\ell},\ \mathbf{C}\in\mathbb{K}^{\ell\times n\times p}\;\text{ and }\;\mathbf{D}\in\mathbb{K}^{m\times n\times p},
\]
be given such that $\ell$ $<$ min$\left(m,n,p\right)$ and 
\[
\mathbf{D}=\mbox{Prod}\left(\mathbf{A},\mathbf{B},\mathbf{C}\right).
\]
The set of hypermatrices
\[
\left\{ \mathbf{X}\in\mathbb{K}^{m\times n\times\ell}\,\text{ such that }\,\mathbf{D}=\mbox{Prod}\left(\mathbf{A},\mathbf{X},\mathbf{C}\right)\right\} 
\]
is determined by solving for the matrix depth slices of $\mathbf{X}$
in the general linear system 
\[
\left\{ \sum_{0\le t<\ell}\mbox{diag}\left(\mathbf{A}\left[:,t,k\right]\right)\cdot\mbox{Mat}\left(\mathbf{X}\left[:,:,t\right]\right)\cdot\mbox{diag}\left(\mathbf{C}\left[t,:,k\right]\right)=\mbox{Mat}\left(\mathbf{D}\left[:,:,k\right]\right)\right\} _{0\le k<p}.
\]
Properties of such general linear constraints motivate our definition
and subsequent investigations of matrix left-right linear diagonal
dependencies. A given multiset of matrices 
\[
\left\{ \mathbf{M}_{t}\right\} _{0\le t<p}\subset\mathbb{K}^{m\times n}
\]
are \emph{left\nobreakdash-right diagonally independent} over the
skew field $\mathbb{K}$ if the only diagonal matrices 
\[
\left\{ \mbox{diag}\left(\mathbf{x}_{t}\right)\right\} _{0\le t<p}\subset\mathbb{K}^{m\times m}\:\text{ and }\:\left\{ \mbox{diag}\left(\mathbf{y}_{t}\right)\right\} _{0\le t<p}\subset\mathbb{K}^{n\times n}
\]
for which the equality
\[
\mathbf{0}_{m\times n}=\sum_{0\le t<p}\mbox{diag}\left(\mathbf{x}_{t}\right)\cdot\mathbf{M}_{t}\cdot\mbox{diag}\left(\mathbf{y}_{t}\right)\text{,}
\]
holds, must be such that
\[
\mathbf{0}_{m\times n}=\mbox{diag}\left(\mathbf{x}_{t}\right)\cdot\mathbf{M}_{t}\cdot\mbox{diag}\left(\mathbf{y}_{t}\right),\quad\forall\:t\in\left\{ 0,1,\cdots,p-1\right\} .
\]
Note that over $\mathbb{C}$, when each matrix in the set $\left\{ \mathbf{M}_{t}\right\} _{0\le t<p}$
has a unique non-zero entry, the left-right diagonal dependence reduces
to the usal notion linear dependence. For such matrices the tight
upper bound for the maximum number of diagonally independent matrices
is $m\cdot n$. Otherwise, when every column ( or respectively every
row ) of each of the matrices in the set $\left\{ \mathbf{M}_{t}\right\} _{0\le t<p}$
is a non-zero, the trivial upper bound on the maximum number of diagonally
independent matrices reduces to min$\left(m,n\right)$, for it suffices
to consider left (or respectively right ) diagonal linear combination
constraints of the form
\[
\mathbf{0}_{m\times n}=\sum_{0\le t<p}\mbox{diag}\left(\mathbf{x}_{t}\right)\cdot\mathbf{M}_{t}\cdot\mathbf{I}_{n},
\]
\[
\text{or}
\]
\[
\mathbf{0}_{m\times n}=\sum_{0\le t<p}\mathbf{I}_{m}\cdot\mathbf{M}_{t}\cdot\mbox{diag}\left(\mathbf{y}_{t}\right).
\]
\\
\\
\textbf{Theorem 3} : Let $\mathbf{H}\in\mathbb{C}^{m\times n\times p}$
be such that 
\[
\text{rank}\left(\mathbf{H}\right)=\ell,
\]
where $0<\ell<$ min$\left(m,n,p\right)$. Let the depth matrix slices
of $\mathbf{H}$ be
\[
\left\{ \mbox{Mat}\left(\mathbf{H}\left[:,:,k\right]\right)\in\mathbb{C}^{m\times n}\right\} _{0\le k<p}
\]
then there exist diagonal matrices $\left\{ \mbox{diag}\left(\mathbf{x}_{it}\right)\right\} _{i,t}\subset\mathbb{C}^{m\times m}$
as well as diagonal matrices $\left\{ \mbox{diag}\left(\mathbf{y}_{tj}\right)\right\} _{t,j}\subset\mathbb{C}^{n\times n}$
such that
\[
\mathbf{0}_{m\times n}=\sum_{0\le t\le\ell}\left(\sum_{0\le i\le\ell-t}\mbox{diag}\left(\mathbf{x}_{i\,t}\right)\cdot\mathbf{M}_{t}\cdot\mbox{diag}\left(\mathbf{y}_{t\,i}\right)\right).
\]
for some choice of a matrix subset $\left\{ \mathbf{M}_{t}\,:\,t\in\left\{ 0,\cdots,\ell\right\} \right\} \subset\left\{ \mbox{Mat}\left(\mathbf{H}\left[:,:,k\right]\right)\,:\,t\in\left\{ 0,\cdots,p-1\right\} \right\} $\\
\\
\emph{Proof} : By the premise that rank$\left(\mathbf{H}\right)=\ell$,
it follows that there is a conformable triple 
\[
\mathbf{U}\in\mathbb{C}^{m\times\ell\times p},\ \mathbf{V}\in\mathbb{C}^{m\times n\times\ell}\ \text{ and }\ \mathbf{W}\in\mathbb{C}^{\ell\times n\times p}
\]
such that 
\[
\mathbf{H}=\mbox{Prod}\left(\mathbf{U},\mathbf{V},\mathbf{W}\right).
\]
Consequently, 
\[
\left\{ \sum_{0\le t<\ell}\mbox{diag}\left(\mathbf{U}\left[:,t,k\right]\right)\cdot\mbox{Mat}\left(\mathbf{V}\left[:,:,t\right]\right)\cdot\mbox{diag}\left(\mathbf{W}\left[t,:,k\right]\right)=\mbox{Mat}\left(\mathbf{\mathbf{H}}\left[:,:,k\right]\right)\right\} _{0\le k<p}.
\]
For fix $\mathbf{U}$ and $\mathbf{W}$ consider the equation in $\mathbf{X}$
of size $m\times n\times\ell$ given by 
\[
\mathbf{H}=\mbox{Prod}\left(\mathbf{U},\mathbf{X},\mathbf{W}\right).
\]
In this form, the equation expresses a general linear system of equation
whose coefficient hypermatrices are $\mathbf{A}\in\mathbb{C}^{1\times\ell\times p}$,
and $\mathbf{B}\in\mathbb{C}^{\ell\times1\times p}$, in an unknown
vector $\mathbf{x}$ of size $1\times1\times\ell$, and the corresponding
right-hand side vector $\mathbf{c}$ of size $1\times1\times p$.
Using the BM product we re-express the corresponding general linear
system of equations as follows 
\begin{equation}
\mbox{Prod}\left(\mathbf{A},\mathbf{x},\mathbf{B}\right)=\mathbf{c}\Leftrightarrow\left\{ \sum_{0\le t<\ell}\mathbf{A}\left[0,t,k\right]\cdot\mathbf{x}\left[0,0,t\right]\cdot\mathbf{B}\left[t,0,k\right]\,=\,\mathbf{c}\left[0,0,k\right]\right\} _{0\le k<p}.\label{Diagonal dependence constraints}
\end{equation}
The entries of the coefficient hypermatrices $\mathbf{A}$ and $\mathbf{B}$
in Eq. (\ref{Diagonal dependence constraints}) are each diagonal
matrices of respective size $m\times m$ and $n\times n$ 
\[
\begin{cases}
\begin{array}{ccc}
\mathbf{A}\left[0,t,k\right] & = & \mbox{diag}\left(\mathbf{U}\left[:,t,k\right]\right)\\
\mathbf{B}\left[t,0,k\right] & = & \mbox{diag}\left(\mathbf{W}\left[t,:,k\right]\right)
\end{array}, & \forall\,\begin{array}{c}
0\le t<\ell\\
0\le k<p
\end{array}\end{cases}.
\]
The entries of $\mathbf{x}$ are unknown $m\times n$ matrices and
the entries of $\mathbf{c}$ correspond to $m\times n$ depth matrix
slices of $\mathbf{H}$ specified by 
\[
\forall\,0\le k<3,\quad\mathbf{c}\left[0,0,k\right]=\mbox{Mat}\left(\mathbf{\mathbf{H}}\left[:,:,k\right]\right).
\]

Since $\mathbf{H}$ has rank exactly $\ell$, for every $k\in\left\{ 0,\cdots,p-1\right\} $
there exists $t\in\left\{ 0,\cdots,\ell-1\right\} $ such that 
\[
\mathbf{0}_{m\times n}\ne\mathbf{A}\left[0,t,k\right]\cdot\mathbf{x}\left[0,0,t\right]\cdot\mathbf{B}\left[t,0,k\right].
\]

We perform on the system 
\[
\left\{ \sum_{0\le t<\ell}\mathbf{A}\left[0,t,k\right]\cdot\mathbf{x}\left[0,0,t\right]\cdot\mathbf{B}\left[t,0,k\right]\,=\,\mathbf{c}\left[0,0,k\right]\right\} _{0\le k<p},
\]
a minor variant of the Gaussian elimination procedure which avoids
division. The procedure is best illustrated by describing the first
round of elimination for a generic system. The first sequence of row
operation is described by 
\begin{equation}
\left(-1\right)\mathbf{A}\left[0,0,k\right]\mbox{R}_{0}\mathbf{B}\left[0,0,k\right]+\mathbf{A}\left[0,0,0\right]\mbox{R}_{k}\mathbf{B}\left[0,0,0\right]\rightarrow\mbox{R}_{k},\quad\forall\,0<k<p.\label{Row operations}
\end{equation}
Although, the very first constraint denoted R$_{0}$ is left unchanged
by the first sequence of row operations, we rewrite R$_{0}$ for consistency
of the unknown variables as follows
\[
\text{R}_{0}\,:\:\sum_{0\le t<\ell}\left(\mathbf{I}_{2}\otimes\mathbf{A}\left[0,t,0\right]\right)\left(\mathbf{I}_{2}\otimes\mathbf{x}\left[0,0,t\right]\right)\left(\mathbf{I}_{2}\otimes\mathbf{B}\left[t,0,0\right]\right)=\left(\mathbf{I}_{2}\otimes\mathbf{c}\left[0,0,0\right]\right),
\]
Following the first sequence of row operations the variable $\mathbf{x}\left[0,0,0\right]$
is eliminated from all but the constraints R$_{0}$ via Eq. (\ref{Row operations})
and yields new constraints of the form
\[
\left\{ \sum_{0<t<\ell}-\mathbf{A}\left[0,0,k\right]\mathbf{A}\left[0,t,0\right]\mathbf{x}\left[0,0,t\right]\mathbf{B}\left[t,0,0\right]\mathbf{B}\left[0,0,k\right]+\mathbf{A}\left[0,0,0\right]\mathbf{A}\left[0,t,k\right]\mathbf{x}\left[0,0,t\right]\mathbf{B}\left[t,0,k\right]\mathbf{B}\left[0,0,0\right]\right.
\]
\[
=
\]
\[
\left.\left(-1\right)\mathbf{A}\left[0,0,k\right]\mathbf{c}\left[0,0,0\right]\mathbf{B}\left[0,0,k\right]+\mathbf{A}\left[0,0,0\right]\mathbf{c}\left[0,0,k\right]\mathbf{B}\left[0,0,0\right]\right\} _{0<k<p}
\]
Equivalently re-expressed as 
\[
\left\{ \sum_{0<t<\ell}\left[\left(-\mathbf{A}\left[0,0,k\right]\mathbf{A}\left[0,t,0\right]\right)\oplus\left(\mathbf{A}\left[0,0,0\right]\mathbf{A}\left[0,t,k\right]\right)\right]\left(\mathbf{I}_{2}\otimes\mathbf{x}\left[0,0,t\right]\right)\left[\left(\mathbf{B}\left[t,0,0\right]\mathbf{B}\left[0,0,k\right]\right)\oplus\left(\mathbf{B}\left[t,0,k\right]\mathbf{B}\left[0,0,0\right]\right)\right]\right.
\]
\[
=
\]
\[
\left.\left(\mathbf{U}_{k}-\mathbf{A}\left[0,0,k\right]\mathbf{c}\left[0,0,0\right]\mathbf{B}\left[0,0,k\right]\right)\oplus\left(\mathbf{A}\left[0,0,0\right]\mathbf{c}\left[0,0,k\right]\mathbf{B}\left[0,0,0\right]-\mathbf{U}_{k}\right)\right\} _{0<k<p},
\]
Note that each $\mathbf{U}_{k}$ denotes an $m\times n$ free variable
matrix. The free variables arise from the fact that given an equality
of the form $y_{0}+y_{1}=a+b$, where $y_{0}$ and $y_{1}$ are unknowns
and $a$ and $b$ are known constants, the equality $y_{0}+y_{1}=a+b$
is equivalent to the assertion $y_{1}=a+t$ in conjunction with the
assertion that $y_{1}=b-t$ for any choice of the free parameter $t$.
The argument is the same for the $\mathbf{U}_{k}$s. The procedure
is thus repeated until the constraints are put in Row Echelon Form
(REF for short). Although the system is non-commutative the procedure
avoids division because conformable diagonal matrices commute.

Having described the proposed variant of the Gaussian elimination
procedure, we discuss properties of the obtained REF. Throughout the
procedure, the right-hand sides are always left\nobreakdash-right
diagonal linear combinations of the previously obtained right hand
side entries. Since there are more constraints then variables (because
$\ell<p$) the procedure must result in at least for $\left(p-\ell\right)$
constraints in the REF whose left hand side express identically zero
functions of the variable entries of $\mathbf{x}$. However the right
hand side need not be an identically zero function of the the free
variables. Using the fact 
\[
\left\{ \mathbf{x}\left[0,0,t\right]=\mbox{Mat}\left(\mathbf{V}\left[:,:,t\right]\right)\right\} _{0\le t<\ell}
\]
are solution to the system it follows that there must exist an assignment
to the free variables which annihilates the right hand side of these
$\left(p-\ell\right)$ constraints.\\
The right hand side expresses in the generic case a left-right diagonal
linear combination of the form
\[
\mathbf{0}_{m\times n}=\sum_{0\le t\le\ell}\left(\sum_{0\le i\le\ell-t}\mbox{diag}\left(\mathbf{x}_{i\,t}\right)\cdot\mathbf{M}_{t}\cdot\mbox{diag}\left(\mathbf{y}_{t\,i}\right)\right),
\]
Thus completing the proof.\\
\\
\\
Note that the variant of Gaussian elimination described in the proof
of Thm. 3 does not fully express solutions to the general linear system,
for two reasons. The first reason is that subsequent row operations
introduce constraints on the previously chosen free variables which
must subsequently solved in order to fully express the solutions to
the system. The second reason result from the fact that diagonal left
and right coefficient matrices of unknowns obtained throughout the
procedure need not be invertible. Subsequent results presented in
this note attempt to address these limitations.

\subsection{Towards hypermatrix rank revealing factorizations}

We emphasize the analogy with the matrix case by first recalling the
property of Gaussian elimination which is central to the LU rank revealing
factorization.\\
\\
\\
\textbf{Theorem 4} : Let $\mathbf{X}\in\mathbb{C}^{m\times\ell}$
and $\mathbf{Y}\in\mathbb{C}^{\ell\times n}$ be matrices such that
for some scalars $\left\{ u_{t}\right\} _{0\le t\ne\tau<\ell}\subset\mathbb{C}$
there exist an index $\tau$ such that, 
\[
\mbox{Prod}_{\boldsymbol{\Delta}^{(\tau)}}\left(\mathbf{X},\mathbf{Y}\right)=\mbox{Prod}\left(\mathbf{X}\left[:,\tau\right],\,\sum_{0\le t\ne\tau<\ell}u_{t}\,\mathbf{Y}\left[t,:\right]\right),
\]
then $\mbox{Prod}\left(\mathbf{X},\mathbf{Y}\right)$ has rank at
most $\left(\ell-1\right)$. \\
\\
\emph{Proof }: Our proof assumes without loss of generality that $\tau=\left(\ell-1\right)$.
The argument is exactly the same for any other choice of $0\le\tau<\ell$.
The product of $\mathbf{X}\in\mathbb{C}^{m\times\ell}$ and $\mathbf{Y}\in\mathbb{C}^{\ell\times n}$
is a sum of $\ell$ outer products. In particular, the outer product
of the last column of $\mathbf{X}$ with the last row of $\mathbf{Y}$
is given by 
\[
\mbox{Prod}_{\boldsymbol{\Delta}^{(\ell-1)}}\left(\mathbf{X},\mathbf{Y}\right)=\mbox{Prod}\left(\mathbf{X}\left[:,\ell-1\right],\mathbf{Y}\left[\ell-1,:\right]\right).
\]
By our hypothesis the row space of $\mathbf{Y}$ has dimension at
most $\left(\ell-1\right)$, and more explicitly we have
\[
\mbox{Prod}_{\boldsymbol{\Delta}^{(\ell-1)}}\left(\mathbf{X},\mathbf{Y}\right)=\mbox{Prod}\left(\mathbf{X}\left[:,\ell-1\right],\,\sum_{0\le t<\ell-1}u_{t}\,\mathbf{Y}\left[t,:\right]\right).
\]
Given the assumption, we express Prod$\left(\mathbf{X},\mathbf{Y}\right)$
as a smaller sum of outer products as follows:
\[
\mbox{Prod}_{\boldsymbol{\Delta}^{(\ell-1)}}\left(\mathbf{X},\mathbf{Y}\right)=\sum_{0\le t<\ell-1}\mbox{Prod}\left(u_{t}\,\mathbf{X}\left[:,\ell-1\right],\,\mathbf{Y}\left[t,:\right]\right),
\]
\[
\implies\mbox{Prod}\left(\mathbf{X},\mathbf{Y}\right)=\left(\sum_{0\le t<\ell-1}\mbox{Prod}_{\boldsymbol{\Delta}^{(t)}}\left(\mathbf{X},\mathbf{Y}\right)\right)+\mbox{Prod}\left(\mathbf{X}\left[:,\ell-1\right],\,\sum_{0\le t<\ell-1}u_{t}\,\mathbf{Y}\left[t,:\right]\right),
\]
and hence
\begin{equation}
\mbox{Prod}\left(\mathbf{X},\mathbf{Y}\right)=\sum_{0\le t<\ell-1}\mbox{Prod}\left(\left(\mathbf{X}\left[:,t\right]+u_{t}\,\mathbf{X}\left[:,\ell-1\right]\right),\,\mathbf{Y}\left[t,:\right]\right).\label{Column linear combination}
\end{equation}
\ref{Column linear combination} expresses the well known elementary
column linear combination operations
\[
\mathbf{X}\left[:,t\right]+u_{t}\,\mathbf{X}\left[:,\ell-1\right]\rightarrow\mathbf{X}\left[:,t\right],\quad\forall\:0\le t<\ell-1.
\]
We now extend the argument above to third order hypermatrices.\\
\\
\\
\textbf{Theorem 5} : Let the conformable triple 
\[
\mathbf{X}\in\mathbb{C}^{m\times\ell\times p},\ \mathbf{Y}\in\mathbb{C}^{m\times n\times\ell}\;\text{ and }\;\mathbf{Z}\in\mathbb{C}^{\ell\times n\times p}
\]
be such that for some index $\tau$ the following holds: 
\[
\forall\:0\le k<\ell,\quad\mbox{Prod}_{\boldsymbol{\Delta}^{(\tau)}}\left(\mathbf{X},\mathbf{Y},\mathbf{Z}\right)\left[:,:,k\right]=
\]
\[
\sum_{0\le t\ne\tau<\ell}\left[\mbox{diag}\left(\mathbf{u}_{t}\right)\cdot\mbox{diag}\left(\mathbf{X}\left[:,\tau,k\right]\right)\cdot\mathbf{Y}\left[:,:,t\right]\cdot\mbox{diag}\left(\mathbf{Z}\left[\tau,:,k\right]\right)\cdot\mbox{diag}\left(\mathbf{v}_{t}\right)+\right.
\]
\[
\mbox{diag}\left(\mathbf{u}_{t}\right)\cdot\mbox{diag}\left(\mathbf{X}\left[:,\tau,k\right]\right)\cdot\mathbf{Y}\left[:,:,t\right]\cdot\mbox{diag}\left(\mathbf{Z}\left[t,:,k\right]\right)+
\]
\[
\left.\mbox{diag}\left(\mathbf{X}\left[:,t,k\right]\right)\cdot\mathbf{Y}\left[:,:,t\right]\cdot\mbox{diag}\left(\mathbf{Z}\left[\tau,:,k\right]\right)\cdot\mbox{diag}\left(\mathbf{v}_{t}\right)\right].
\]
then $\mbox{Prod}\left(\mathbf{X},\mathbf{Y},\mathbf{Z}\right)$ has
rank at most $\left(\ell-1\right)$. \\
\\
\\
\emph{Proof }: Our proof assumes without loss of generality that $\tau=\left(\ell-1\right)$.
The argument is exactly the same for any choice $0\le\tau<\ell$.
The product of the conformable triple 
\[
\mathbf{X}\in\mathbb{C}^{m\times\ell\times p},\ \mathbf{Y}\in\mathbb{C}^{m\times n\times\ell}\;\text{ and }\;\mathbf{Z}\in\mathbb{C}^{\ell\times n\times p}
\]
expresses a sum of $\ell$ outer products. In particular, consider
the outer product of the last column slice of $\mathbf{X}$ with the
last depth slice of $\mathbf{Y}$ and the last row slice of $\mathbf{Z}$
given by
\[
\mbox{Prod}_{\boldsymbol{\Delta}^{(\ell-1)}}\left(\mathbf{X},\mathbf{Y},\mathbf{Z}\right)=\mbox{Prod}\left(\mathbf{X}\left[:,\ell-1,:\right],\mathbf{Y}\left[:,:,\ell-1\right],\mathbf{Z}\left[\ell-1,:,:\right]\right).
\]
By analogy to the matrix case, an assumption is made on this outer
product 
\[
\mbox{Prod}\left(\mathbf{X}\left[:,\ell-1,:\right],\mathbf{Y}\left[:,:,\ell-1\right],\mathbf{Z}\left[\ell-1,:,:\right]\right).
\]
The explicit assumption is that 
\[
\forall\:0\le k<\ell,\quad\mbox{Prod}_{\boldsymbol{\Delta}^{(\ell-1)}}\left(\mathbf{X},\mathbf{Y},\mathbf{Z}\right)\left[:,:,k\right]=
\]
\[
\sum_{0\le t<\ell-1}\left[\mbox{diag}\left(\mathbf{u}_{t}\right)\cdot\mbox{diag}\left(\mathbf{X}\left[:,\ell-1,k\right]\right)\cdot\mathbf{Y}\left[:,:,t\right]\cdot\mbox{diag}\left(\mathbf{Z}\left[\ell-1,:,k\right]\right)\cdot\mbox{diag}\left(\mathbf{v}_{t}\right)+\right.
\]
\[
\mbox{diag}\left(\mathbf{u}_{t}\right)\cdot\mbox{diag}\left(\mathbf{X}\left[:,\ell-1,k\right]\right)\cdot\mathbf{Y}\left[:,:,t\right]\cdot\mbox{diag}\left(\mathbf{Z}\left[t,:,k\right]\right)+
\]
\begin{equation}
\left.\mbox{diag}\left(\mathbf{X}\left[:,t,k\right]\right)\cdot\mathbf{Y}\left[:,:,t\right]\cdot\mbox{diag}\left(\mathbf{Z}\left[\ell-1,:,k\right]\right)\cdot\mbox{diag}\left(\mathbf{v}_{t}\right)\right].\label{dependence assumption}
\end{equation}
Given our assumption we express Prod$\left(\mathbf{X},\mathbf{Y},\mathbf{Z}\right)$
as a smaller sum of outer products as follows:

\[
\forall\:0\le k<p,\quad\mbox{Prod}\left(\mathbf{X},\,\mathbf{Y},\,\mathbf{Z}\right)\left[:,:,k\right]=\sum_{0\le t<\ell-1}\mbox{diag}\left(\mathbf{X}\left[:,t,k\right]\right)\cdot\mathbf{Y}\left[:,:,t\right]\cdot\mbox{diag}\left(\mathbf{Z}\left[t,:,k\right]\right)+
\]
\[
\sum_{0\le t<\ell-1}\left[\mbox{diag}\left(\mathbf{u}_{t}\right)\mbox{diag}\left(\mathbf{X}\left[:,\ell-1,k\right]\right)\mathbf{Y}\left[:,:,t\right]\mbox{diag}\left(\mathbf{Z}\left[\ell-1,:,k\right]\right)\mbox{diag}\left(\mathbf{v}_{t}\right)+\right.
\]
\[
\mbox{diag}\left(\mathbf{u}_{t}\right)\mbox{diag}\left(\mathbf{X}\left[:,\ell-1,k\right]\right)\mathbf{Y}\left[:,:,t\right]\mbox{diag}\left(\mathbf{Z}\left[t,:,k\right]\right)+
\]
\[
\left.\mbox{diag}\left(\mathbf{X}\left[:,t,k\right]\right)\mathbf{Y}\left[:,:,t\right]\mbox{diag}\left(\mathbf{Z}\left[\ell-1,:,k\right]\right)\mbox{diag}\left(\mathbf{v}_{t}\right)\right],
\]
hence
\[
\forall\:0\le k<p,\quad\mbox{Prod}\left(\mathbf{X},\mathbf{Y},\mathbf{Z}\right)\left[:,:,k\right]=
\]
\[
\sum_{0\le t<\ell-1}\left(\mbox{diag}\left(\mathbf{u}_{t}\right)\mbox{diag}\left(\mathbf{X}\left[:,\ell-1,k\right]\right)+\mbox{diag}\left(\mathbf{X}\left[:,t,k\right]\right)\right)\mathbf{Y}\left[:,:,t\right]\left(\mbox{diag}\left(\mathbf{Z}\left[t,:,k\right]\right)+\mbox{diag}\left(\mathbf{Z}\left[\ell-1,:,k\right]\right)\mbox{diag}\left(\mathbf{v}_{t}\right)\right).
\]
Thus completing our proof.\\
\\
The condition described in Thm. (4) for the matrix case and Thm. (5)
for the hypermatrix case are sufficient but clearly not necessary.
To derive from Thm. (4), a necessary condition for the rank of the
matrix $\mathbf{Y}\in\mathbb{C}^{\ell\times n}$ to be less than $\ell$,
it suffices to the set $\mathbf{X}$ to the identity matrix $\mathbf{I}_{\ell}$.
Similarly, we derive from Thm. (5) a necessary condition for the rank
of the hypermatrix $\mathbf{Y}\in\mathbb{C}^{m\times\ell\times p}$
to be rank less then $\ell$, it suffices to set the hypermatrices
$\mathbf{X}$ and $\mathbf{Y}$ to be the identity pair $\mathbf{J}_{0}\in\mathbb{C}^{m\times p\times p}$
and $\mathbf{J}_{1}\in\mathbb{C}^{p\times n\times p}$ whose entries
are specified by 
\[
\mathbf{J}_{0}\left[i,t,k\right]=\begin{cases}
\begin{array}{cc}
1 & \mbox{ if }0\le t=k<p\\
0 & \mbox{otherwise}
\end{array} & \forall\,\begin{cases}
\begin{array}{c}
0\le i<m\\
0\le t<p\\
0\le k<p
\end{array}\end{cases}\end{cases},
\]
\[
\mathbf{J}_{1}\left[t,j,k\right]=\begin{cases}
\begin{array}{cc}
1 & \mbox{ if }0\le t=k<p\\
0 & \mbox{otherwise}
\end{array} & \forall\,\begin{cases}
\begin{array}{c}
0\le t<p\\
0\le j<n\\
0\le k<p
\end{array}\end{cases}\end{cases},
\]
For which we recall the defining property being that 
\[
\forall\:\mathbf{A}\in\mathbb{C}^{m\times n\times p},\quad\mathbf{A}=\mbox{Prod}\left(\mathbf{J}_{0},\mathbf{A},\mathbf{J}_{1}\right).
\]
where we assume that $n=$ min$\left(m,\,n,\,p\right)$. And hence,
the condition in Eq. (\ref{dependence assumption}) therefore reduces
to left and right diagonal dependence of the depth slices.\\
\\
The rank revealing factorizations of Prod$\left(\mathbf{X},\mathbf{Y},\mathbf{Z}\right)$
is obtained by repeatedly solving for entries of $\left\{ \mbox{diag}\left(\mathbf{u}_{t}\right),\,\mbox{diag}\left(\mathbf{v}_{t}\right)\right\} _{0\le t<\ell-1}$
in Eq. (\ref{dependence assumption}). At each iteration the assignment
of $\left\{ \mbox{diag}\left(\mathbf{u}_{t}\right),\,\mbox{diag}\left(\mathbf{v}_{t}\right)\right\} _{0\le t<\ell-1}$
effectively reduces by one the number of outer product summands. The
corresponding elementary slice operation is thus prescribed by
\[
\begin{cases}
\begin{array}{ccc}
\left(\mbox{diag}\left(\mathbf{u}_{t}\right)\mathbf{X}\left[:,\ell-1,k\right]+\mathbf{X}\left[:,t,k\right]\right) & \rightarrow & \mathbf{X}\left[:,t,k\right]\\
\left(\mathbf{Z}\left[t,:,k\right]+\mathbf{Z}\left[\ell-1,:,k\right]\mbox{diag}\left(\mathbf{v}_{t}\right)\right) & \rightarrow & \mathbf{Z}\left[t,:,k\right]
\end{array} & \,\forall\:0\le t<\ell-1\end{cases}.
\]
Unfortunately the constraints in Eq. (\ref{dependence assumption})
are non-linear constraints.\\
\\
\\
\textbf{Theorem 6} : If a generic $\mathbf{B}\in\mathbb{C}^{m\times n\times\left(r+1\right)}$,
such that $\mathbf{B}\left[:,:,0\right]\in\left(\mathbb{C}-\left\{ 0\right\} \right)^{m\times n\times1}$
and $r<\min\left(m,n\right)$, and then 
\[
\mbox{rank}\left(\mathbf{B}\right)=\left(r+1\right)\implies\left(m+n\right)\cdot\left(r-1\right)<\left(m-1\right)\cdot\left(n-1\right).
\]
\\
\emph{Proof} : It suffices to focus on the case 
\[
\mbox{rank}\left(\mathbf{B}\right)=r.
\]
It follows from Theorem 5 that the outer product reduction criterion
for the product of a conformable triple of hypermatrices 
\[
\mathbf{A}\in\mathbb{C}^{m\times\left(r+1\right)\times\left(r+1\right)},\ \mathbf{B}\in\mathbb{C}^{m\times n\times\left(r+1\right)}\ \text{ and }\ \mathbf{C}\in\mathbb{C}^{\left(r+1\right)\times n\times\left(r+1\right)}
\]
is given by 
\[
\forall\:0\le k<n,\quad\mbox{Prod}_{\boldsymbol{\Delta}^{(\tau)}}\left(\mathbf{A},\mathbf{B},\mathbf{C}\right)\left[:,:,k\right]=
\]
\[
\sum_{0\le t\ne\tau<r+1}\left[\mbox{diag}\left(\mathbf{u}_{t}\right)\mbox{diag}\left(\mathbf{A}\left[:,\tau,k\right]\right)\mathbf{B}\left[:,:,t\right]\mbox{diag}\left(\mathbf{C}\left[\tau,:,k\right]\right)\mbox{diag}\left(\mathbf{v}_{t}\right)+\right.
\]
\[
\mbox{diag}\left(\mathbf{u}_{t}\right)\cdot\mbox{diag}\left(\mathbf{A}\left[:,\tau,k\right]\right)\mathbf{B}\left[:,:,t\right]\mbox{diag}\left(\mathbf{C}\left[t,:,k\right]\right)+
\]
\[
\left.\mbox{diag}\left(\mathbf{A}\left[:,t,k\right]\right)\mathbf{B}\left[:,:,t\right]\mbox{diag}\left(\mathbf{C}\left[\tau,:,k\right]\right)\mbox{diag}\left(\mathbf{v}_{t}\right)\right].
\]
where $\mathbf{A}=\mathbf{J}_{0}\in\mathbb{C}^{m\times\left(r+1\right)\times\left(r+1\right)}$
and $\mathbf{C}=\mathbf{J}_{1}\in\mathbb{C}^{\left(r+1\right)\times n\times\left(r+1\right)}$
whose entries are specified by
\[
\mathbf{J}_{0}\left[i,t,k\right]=\begin{cases}
\begin{array}{cc}
1 & \mbox{ if }0\le t=k<r+1\\
0 & \mbox{otherwise}
\end{array} & \forall\,\begin{cases}
\begin{array}{c}
0\le i<m\\
0\le t<r+1\\
0\le k<r+1
\end{array}\end{cases}\end{cases},
\]
\[
\mathbf{J}_{1}\left[t,j,k\right]=\begin{cases}
\begin{array}{cc}
1 & \mbox{ if }0\le t=k<r+1\\
0 & \mbox{otherwise}
\end{array} & \forall\,\begin{cases}
\begin{array}{c}
0\le t<r+1\\
0\le j<n\\
0\le k<r+1
\end{array}\end{cases}\end{cases},
\]
we have 
\[
\mbox{Prod}\left(\mathbf{A},\mathbf{B},\mathbf{C}\right)=\mbox{Prod}\left(\mathbf{J}_{0},\mathbf{B},\mathbf{J}_{1}\right)=\mathbf{B},
\]
and the equality 
\[
\forall\:0\le k<n,\quad\mbox{Prod}_{\boldsymbol{\Delta}^{(\tau)}}\left(\mathbf{A},\mathbf{B},\mathbf{C}\right)\left[:,:,k\right]=
\]
\[
\sum_{0\le t\ne\tau<r+1}\left[\mbox{diag}\left(\mathbf{u}_{t}\right)\mbox{diag}\left(\mathbf{A}\left[:,\tau,k\right]\right)\mathbf{B}\left[:,:,t\right]\mbox{diag}\left(\mathbf{C}\left[\tau,:,k\right]\right)\mbox{diag}\left(\mathbf{v}_{t}\right)+\right.
\]
\[
\mbox{diag}\left(\mathbf{u}_{t}\right)\cdot\mbox{diag}\left(\mathbf{A}\left[:,\tau,k\right]\right)\mathbf{B}\left[:,:,t\right]\mbox{diag}\left(\mathbf{C}\left[t,:,k\right]\right)+
\]
\[
\left.\mbox{diag}\left(\mathbf{A}\left[:,t,k\right]\right)\mathbf{B}\left[:,:,t\right]\mbox{diag}\left(\mathbf{C}\left[\tau,:,k\right]\right)\mbox{diag}\left(\mathbf{v}_{t}\right)\right]
\]
simplifies to the following diagonal dependence relation
\[
\mathbf{B}\left[:,:,\tau\right]=\sum_{0\le t\ne\tau<r+1}\mbox{diag}\left(\mathbf{U}\left[:,t\right]\right)\mathbf{B}\left[:,:,t\right]\mbox{diag}\left(\mathbf{V}\left[t,:\right]\right).
\]
which we re-express in terms of columns and rows of $\mathbf{U}\in\mathbb{C}^{m\times p}$
and $\mathbf{V}^{p\times n}$ respectively as 
\[
\mathbf{0}_{n\times n}=\sum_{0\le t<r+1}\mbox{diag}\left(\mathbf{U}\left[:,t\right]\right)\mathbf{B}\left[:,:,t\right]\mbox{diag}\left(\mathbf{V}\left[t,:\right]\right).
\]
\begin{equation}
\mbox{Prod}\left(\mathbf{X}^{\prime},\mathbf{B},\mathbf{Y}^{\prime}\right)=\mathbf{0}_{m\times n\times1},\label{Dependence Constraints}
\end{equation}
where $\mathbf{X}^{\prime}\in\mathbb{C}^{m\times\left(r+1\right)\times1}$
and $\mathbf{Y}^{\prime}\in\mathbb{C}^{\left(r+1\right)\times n\times1}$,
Let $\tau=0$ then we have 
\[
0\ne\det\left\{ \mbox{diag}\left(\mathbf{X}^{\prime}\left[:,r,0\right]\right)\right\} \det\left\{ \mbox{diag}\left(\mathbf{Y}^{\prime}\left[r,:,0\right]\right)\right\} .
\]
we rewrite the constraints in Eq. (\ref{Dependence Constraints})
as
\[
\sum_{0\le t<r}\left(\frac{\mathbf{X}^{\prime}\left[i,t,0\right]}{\mathbf{X}^{\prime}\left[i,r,0\right]}\right)\mathbf{B}\left[i,j,t\right]\left(\frac{\mathbf{Y}^{\prime}\left[t,j,0\right]}{\mathbf{Y}^{\prime}\left[r,j,0\right]}\right)=-\mathbf{B}\left[i,j,r\right],\;\forall\:\begin{cases}
\begin{array}{c}
0\le i<m\\
0\le j<n
\end{array}\end{cases}.
\]
Hence the constraints in Eq. (\ref{Dependence Constraints}) can be
re-expressed in terms of smaller hypermatrices $\mathbf{X}\in\mathbb{C}^{m\times r\times1}$
and $\mathbf{Y}\in\mathbb{C}^{r\times n\times1}$ as 
\[
\mbox{Prod}\left(\mathbf{X},\mathbf{B}\left[:,:,:r-1\right],\mathbf{Y}\right)=-\mathbf{B}\left[:,:,r\right]\Leftrightarrow\left\{ \sum_{0\le t<r}\mathbf{X}\left[i,t,0\right]\mathbf{B}\left[i,j,t\right]\mathbf{Y}\left[t,j,0\right]=-\mathbf{B}\left[i,j,r\right]\right\} _{\begin{array}{c}
0\le i<m\\
0\le j<n
\end{array}}.
\]
We further cast these constraints as a system of linear equations
of the form
\[
\sum_{0\le n\cdot u+v<m\cdot n}\mathbf{M}\left[n\cdot i+j,n\cdot u+v\right]\left(\mathbf{X}\left[u,0,0\right]\mathbf{Y}\left[0,v,0\right]\right)=
\]
\[
-\left(\mathbf{B}\left[i,j,r\right]+\sum_{0<t<r+1}\mathbf{X}\left[i,t,0\right]\mathbf{B}\left[i,j,t\right]\mathbf{Y}\left[t,j,0\right]\right),\;\forall\:\begin{cases}
\begin{array}{c}
0\le i<m\\
0\le j<n
\end{array}\end{cases},
\]
where 
\[
\mathbf{M}\left[n\cdot i+j,n\cdot u+v\right]=\begin{cases}
\begin{array}{cc}
\mathbf{B}\left[i,j,0\right] & \text{ if }n\cdot i+j=n\cdot u+v\\
0 & \text{otherwise}
\end{array}\end{cases}
\]
\[
\implies\mathbf{X}\left[i,0,0\right]\mathbf{Y}\left[0,j,0\right]=-\left(\frac{\mathbf{B}\left[i,j,r\right]}{\mathbf{B}\left[i,j,0\right]}+\sum_{0<t<r}\mathbf{X}\left[i,t,0\right]\frac{\mathbf{B}\left[i,j,t\right]}{\mathbf{B}\left[i,j,0\right]}\mathbf{Y}\left[t,j,0\right]\right),\;\forall\:\begin{cases}
\begin{array}{c}
0\le i<m\\
0\le j<n
\end{array}\end{cases}.
\]
These algebraic relations effectively eliminate the variables $\left\{ \mathbf{X}\left[u,0,0\right]\right\} _{0\le u<m}$
as well as the variables $\left\{ \mathbf{Y}\left[0,v,0\right]\right\} _{0\le v<n}$
from the system. Note that the monomials $\left\{ \mathbf{X}\left[u,0,0\right]\mathbf{Y}\left[0,v,0\right]\right\} _{\begin{array}{c}
0\le u<m\\
0\le v<n
\end{array}}$ correspond to entries of an $m\times n$ rank one matrix. The problem
is thus reduced to a system of ${m \choose 2}\cdot{n \choose 2}$
constraints in the remaining $\left(m+n\right)\cdot\left(r-1\right)$
variable entries of $\mathbf{X}$ and $\mathbf{Y}$. These constraints
are determinantal constraints of the form
\[
\forall\:\begin{array}{c}
0\le i_{0}<i_{1}<m\\
0\le j_{0}<j_{1}<n
\end{array},\;
\]
 
\begin{equation}
\det\left(\begin{array}{ccc}
\frac{\mathbf{B}\left[i_{0},j_{0},r\right]}{\mathbf{B}\left[i_{0},j_{0},0\right]}+\underset{0<t<r}{\sum}\mathbf{X}\left[i_{0},t,0\right]\frac{\mathbf{B}\left[i_{0},j_{0},t\right]}{\mathbf{B}\left[i_{0},j_{0},0\right]}\mathbf{Y}\left[t,j_{0},0\right] &  & \frac{\mathbf{B}\left[i_{0},j_{1},r\right]}{\mathbf{B}\left[i_{0},j_{1},0\right]}+\underset{0<t<r}{\sum}\mathbf{X}\left[i_{0},t,0\right]\frac{\mathbf{B}\left[i_{0},j_{1},t\right]}{\mathbf{B}\left[i_{0},j_{1},0\right]}\mathbf{Y}\left[t,j_{1},0\right]\\
\\
\frac{\mathbf{B}\left[i_{1},j_{0},r\right]}{\mathbf{B}\left[i_{1},j_{0},0\right]}+\underset{0<t<r}{\sum}\mathbf{X}\left[i_{1},t,0\right]\frac{\mathbf{A}\left[i_{1},j_{0},t\right]}{\mathbf{A}\left[i_{1},j_{0},0\right]}\mathbf{Y}\left[t,j_{0},0\right] &  & \frac{\mathbf{B}\left[i_{1},j_{1},r\right]}{\mathbf{B}\left[i_{1},j_{1},0\right]}+\underset{0<t<r}{\sum}\mathbf{X}\left[i_{1},t,0\right]\frac{\mathbf{B}\left[i_{1},j_{1},t\right]}{\mathbf{B}\left[i_{1},j_{1},0\right]}\mathbf{Y}\left[t,j_{1},0\right]
\end{array}\right)=0.\label{Determinant Constraints}
\end{equation}
The system above is of type 2 as described in \cite{2018arXiv180803743G}
in the $mn$ variables 
\[
\left\{ \frac{\mathbf{B}\left[i,j,r\right]}{\mathbf{B}\left[i,j,0\right]}+\sum_{0<t<r}\mathbf{X}\left[i,t,0\right]\frac{\mathbf{B}\left[i,j,t\right]}{\mathbf{B}\left[i,j,0\right]}\mathbf{Y}\left[t,j,0\right]\right\} _{\begin{array}{c}
0\le i<m\\
0\le j<n
\end{array}}
\]
In particular the exponent matrix of the corresponding REF of the
system has $\left(m-1\right)\left(n-1\right)$ pivot ( accounting
for the degrees of freedom of the rank 1 matrix ) variables, $m+n-1$
free variables and ${m \choose 2}{n \choose 2}-\left(m-1\right)\cdot\left(n-1\right)$
zero rows. By the algebraic independence of the pivot rows it follows
that $\mathbf{A}$ stands a chance to have rank $r+1$ if $\left(m-1\right)\cdot\left(n-1\right)$
exceeds the number of remaining variables. The condition is thus expressed
by the inequality
\[
\left(m+n\right)\cdot\left(r-1\right)<\left(m-1\right)\cdot\left(n-1\right).
\]
\\
\\
\textbf{Theorem 7} : The Bhattacharya-Mesner rank of a generic \footnote{A generic hypermatrix is one whose entries do not satisfy any non-trivial
algebraic relation. In particular, all of its entries are non-zero.} $\mathbf{B}\in\mathbb{C}^{n\times n\times n}$ is at most $2$ if
$n=2$ and at most $\left(n-1\right)$ if $n>2$.\\
\\
\emph{Proof }: It follows from Theorem 5 that the outer product reduction
criterion for the product of a conformable triple of hypermatrices
\[
\mathbf{A}\in\mathbb{C}^{m\times p\times p},\ \mathbf{B}\in\mathbb{C}^{m\times n\times p}\ \text{ and }\ \mathbf{C}\in\mathbb{C}^{p\times n\times p}
\]
where $p=$ min$\left(m,n\right)$ is given by 
\[
\forall\:0\le k<n,\quad\mbox{Prod}_{\boldsymbol{\Delta}^{(\tau)}}\left(\mathbf{A},\mathbf{B},\mathbf{C}\right)\left[:,:,k\right]=
\]
\[
\sum_{0\le t\ne\tau<p}\left[\mbox{diag}\left(\mathbf{u}_{t}\right)\mbox{diag}\left(\mathbf{A}\left[:,\tau,k\right]\right)\mathbf{B}\left[:,:,t\right]\mbox{diag}\left(\mathbf{C}\left[\tau,:,k\right]\right)\mbox{diag}\left(\mathbf{v}_{t}\right)+\right.
\]
\[
\mbox{diag}\left(\mathbf{u}_{t}\right)\cdot\mbox{diag}\left(\mathbf{A}\left[:,\tau,k\right]\right)\mathbf{B}\left[:,:,t\right]\mbox{diag}\left(\mathbf{C}\left[t,:,k\right]\right)+
\]
\[
\left.\mbox{diag}\left(\mathbf{A}\left[:,t,k\right]\right)\mathbf{B}\left[:,:,t\right]\mbox{diag}\left(\mathbf{C}\left[\tau,:,k\right]\right)\mbox{diag}\left(\mathbf{v}_{t}\right)\right].
\]
where $\mathbf{A}=\mathbf{J}_{0}\in\mathbb{C}^{m\times p\times p}$
and $\mathbf{C}=\mathbf{J}_{1}\in\mathbb{C}^{p\times n\times p}$
whose entries are specified by
\[
\mathbf{J}_{0}\left[i,t,k\right]=\begin{cases}
\begin{array}{cc}
1 & \mbox{ if }0\le t=k<p\\
0 & \mbox{otherwise}
\end{array} & \forall\,\begin{cases}
\begin{array}{c}
0\le i<m\\
0\le t<p\\
0\le k<p
\end{array}\end{cases}\end{cases},
\]
\[
\mathbf{J}_{1}\left[t,j,k\right]=\begin{cases}
\begin{array}{cc}
1 & \mbox{ if }0\le t=k<p\\
0 & \mbox{otherwise}
\end{array} & \forall\,\begin{cases}
\begin{array}{c}
0\le t<p\\
0\le j<n\\
0\le k<p
\end{array}\end{cases}\end{cases},
\]
we have 
\[
\mbox{Prod}\left(\mathbf{A},\mathbf{B},\mathbf{C}\right)=\mbox{Prod}\left(\mathbf{J}_{0},\mathbf{B},\mathbf{J}_{1}\right)=\mathbf{B},
\]
and the equality 
\[
\forall\:0\le k<n,\quad\mbox{Prod}_{\boldsymbol{\Delta}^{(\tau)}}\left(\mathbf{A},\mathbf{B},\mathbf{C}\right)\left[:,:,k\right]=
\]
\[
\sum_{0\le t\ne\tau<p}\left[\mbox{diag}\left(\mathbf{u}_{t}\right)\mbox{diag}\left(\mathbf{A}\left[:,\tau,k\right]\right)\mathbf{B}\left[:,:,t\right]\mbox{diag}\left(\mathbf{C}\left[\tau,:,k\right]\right)\mbox{diag}\left(\mathbf{v}_{t}\right)+\right.
\]
\[
\mbox{diag}\left(\mathbf{u}_{t}\right)\cdot\mbox{diag}\left(\mathbf{A}\left[:,\tau,k\right]\right)\mathbf{B}\left[:,:,t\right]\mbox{diag}\left(\mathbf{C}\left[t,:,k\right]\right)+
\]
\[
\left.\mbox{diag}\left(\mathbf{A}\left[:,t,k\right]\right)\mathbf{B}\left[:,:,t\right]\mbox{diag}\left(\mathbf{C}\left[\tau,:,k\right]\right)\mbox{diag}\left(\mathbf{v}_{t}\right)\right]
\]
simplifies to the following diagonal dependence relation
\begin{equation}
\mathbf{B}\left[:,:,\tau\right]=\sum_{0\le t\ne\tau<p}\mbox{diag}\left(\mathbf{U}\left[:,t\right]\right)\mathbf{B}\left[:,:,t\right]\mbox{diag}\left(\mathbf{V}\left[t,:\right]\right).\label{Slice dependence}
\end{equation}
which we re-express in terms of columns and rows of $\mathbf{U}\in\mathbb{C}^{m\times p}$
and $\mathbf{V}^{p\times n}$ respectively as The diagonal dependence
is therefore expressed by constraints of the form 
\[
\mathbf{0}_{m\times n}=\sum_{0\le t<p}\mbox{diag}\left(\mathbf{U}\left[:,t\right]\right)\mathbf{B}\left[:,:,t\right]\mbox{diag}\left(\mathbf{V}\left[t,:\right]\right).
\]
Consequently for all $1\le k<n$ we have
\[
\mathbf{U}\left[:,k\right]\mathbf{V}\left[k,:\right]=\sum_{0\le t\ne k<p}\mbox{diag}\left(\mathbf{U}\left[:,t\right]\right)\left(\mathbf{B}\left[:,:,t\right]\circ\left(\mathbf{B}\left[:,:,k\right]\right)^{\circ^{-1}}\right)\mbox{diag}\left(\mathbf{V}\left[t,:\right]\right)
\]
\[
\implies\mathbf{U}\cdot\mathbf{V}=\sum_{0\le k<p}\left(\sum_{1\le t\ne k<p}\mbox{diag}\left(\mathbf{U}\left[:,t\right]\right)\left(\mathbf{B}\left[:,:,t\right]\circ\left(\mathbf{B}\left[:,:,k\right]\right)^{\circ^{-1}}\right)\mbox{diag}\left(\mathbf{V}\left[t,:\right]\right)\right).
\]
Here $\circ$ denotes the entry-wise product also called the Hadamard
product, and $\circ^{k}$ denotes the entry-wise exponentiation by
$k$ of non-zero entries. The diagonal dependence status of the depth
slices of $\mathbf{B}$ does not change if each depth slice of $\mathbf{B}$
is premultiplied on the right by distinct known invertible diagonal
matrices. As a result, we can assume without loss of generality that
det$\left(\mathbf{V}\right)\ne0$ and $n\le m$, in which we case
we have
\begin{equation}
\mathbf{U}\det\left(\mathbf{V}\right)=\left(\sum_{0\le k<p}\sum_{0\le t\ne k<p}\mbox{diag}\left(\mathbf{U}\left[:,t\right]\right)\left(\mathbf{B}\left[:,:,t\right]\circ\left(\mathbf{B}\left[:,:,k\right]\right)^{\circ^{-1}}\right)\mbox{diag}\left(\mathbf{V}\left[t,:\right]\right)\right)\text{adj}\left(\mathbf{V}\right).\label{e-vector/e-value}
\end{equation}
\[
\implies\mathbf{U}\det\left(\mathbf{V}\right)=\left(\sum_{0\le k<p}\sum_{0\le t\ne k<p}\left(\mathbf{U}\left[:,t\right]\cdot\mathbf{V}\left[t,:\right]\right)\circ\left(\mathbf{B}\left[:,:,t\right]\circ\left(\mathbf{B}\left[:,:,k\right]\right)^{\circ^{-1}}\right)\right)\cdot\text{adj}\left(\mathbf{V}\right).
\]
Let $\mathbf{M}\in\left(\mathbb{C}\left[v_{00},\cdots,v_{n-1\,n-1}\right]\right)^{mp\times mp}$
whose entries are given by 
\[
\mathbf{M}\left[p\cdot i+j,\,p\cdot s+t\right]=\frac{\partial}{\partial\mathbf{U}\left[s,t\right]}\left(\left(\sum_{0\le k<p}\sum_{0\le t\ne k<p}\left(\mathbf{U}\left[:,t\right]\cdot\mathbf{V}\left[t,:\right]\right)\circ\left(\mathbf{B}\left[:,:,t\right]\circ\left(\mathbf{B}\left[:,:,k\right]\right)^{\circ^{-1}}\right)\right)\cdot\text{adj}\left(\mathbf{V}\right)\right)\left[i,j\right]
\]
The constraints in Eq. (\ref{e-vector/e-value}) express an eigenvector/eigenvalue
problem for the $m^{2}\times m^{2}$ matrix $\mathbf{M}$, whose entries
are polynomials in the entries of $\mathbf{V}$. The desired eigenvalue
is det$\left(\mathbf{V}\right)$ while the entries of $\mathbf{U}$
make up the entries of the eigenvector. The constraints are of the
form 
\[
\left\{ \mathbf{U}\left[i,j\right]\det\left(\mathbf{V}\right)=\sum_{\begin{array}{c}
0\le s<m\\
0\le t<p
\end{array}}\mathbf{M}\left[p\cdot i+j,\,p\cdot s+t\right]\mathbf{U}\left[s,t\right]\right\} _{\begin{array}{c}
0\le i<m\\
0\le j<p
\end{array}}.
\]
In order to find the eigenvalues, we solve for the characteristic
polynomial equation 
\[
\det\left\{ \det\left(\mathbf{V}\right)\mathbf{I}_{mp}-\mathbf{M}\right\} =\prod_{0\le i<m}\det\left\{ \det\left(\mathbf{V}\right)\mathbf{I}_{p}-\mathbf{M}\left[i\,p:\left(i+1\right)p,\,i\,p:\left(i+1\right)p\right]\right\} .
\]
For generic $\mathbf{B}$ the polynomial $\det\left\{ \det\left(\mathbf{V}\right)\mathbf{I}_{mp}-\mathbf{M}\right\} \in\mathbb{C}\left[v_{00},\cdots,v_{p-1\,p-1}\right]$
is not an identically constant polynomial, for if this were the case
then the entries of $\mathbf{B}$ would satisfy some non-trivial algebraic
relation. Consequently the polynomial det$\left\{ \det\left(\mathbf{V}\right)\mathbf{I}_{mp}-\mathbf{M}\right\} $
admits roots in $\mathbb{C}$. Having assumed that det$\left(\mathbf{V}\right)\ne0$
and that the entries of $\mathbf{B}$ are generic, the constraints
reduce to a system of equations of the form 
\[
\left\{ 0=\sum_{\sigma\in\text{S}_{n}}f_{j,\sigma}\left(\mathbf{B}\right)\,\text{sgn}\left(\sigma\right)\prod_{0\le i<n}\mathbf{V}\left[i,\sigma\left(i\right)\right]\right\} _{0\le j<m}
\]
where 
\[
\left\{ f_{j,\sigma}\left(\mathbf{B}\right)\,:\,0\le i<m,\:\sigma\in\text{S}_{n}\right\} \subset\mathbb{C}\left[b_{000},\cdots,b_{m-1\,n-1\,p-1}\right].
\]
Each constraint therefore translates into constraints of the form
\[
\left\{ \prod_{0\le i<n}\mathbf{V}\left[i,\sigma\left(i\right)\right]=\sum_{\gamma\in\text{S}_{n}\backslash\left\{ \text{id}\right\} }\frac{c_{j,\gamma}\,\text{sgn}\left(\gamma\right)}{f_{j,\gamma}\left(\mathbf{B}\right)}\exp\left(\frac{2\pi\,\text{lex}\sigma\cdot\text{lex}\gamma}{n!}\sqrt{-1}\right)\right\} _{\sigma\in\text{S}_{n}}
\]
where lex denotes the bijective map from the symmetric group S$_{n}$
to $\left\{ 0,\cdots,n!-1\right\} $ which reflects the lexicographical
ordering of the permutation strings. Note that in the case $n=2$
we have 
\[
0=f_{j,\text{id}}\left(\mathbf{B}\right)v_{00}v_{11}+\left(-1\right)f_{j,1-\text{id}}\left(\mathbf{B}\right)v_{01}v_{10}\implies\left(\begin{array}{c}
v_{00}v_{11}\\
\\
v_{01}v_{10}
\end{array}\right)=\left(\begin{array}{c}
\frac{c_{j,1-\text{id}}}{f_{i,\text{id}}\left(\mathbf{B}\right)}\\
\\
\frac{c_{j,1-\text{id}}}{f_{i,1-\text{id}}\left(\mathbf{B}\right)}
\end{array}\right)
\]
for some free variable $c_{1-\text{id}}$ 
\[
\implies\begin{cases}
\begin{array}{ccc}
v_{00} & = & \frac{c_{j,1-\text{id}}\,s_{j}^{-1}}{f_{i,\text{id}}\left(\mathbf{B}\right)}\\
\\
v_{01} & = & \frac{c_{j,1-\text{id}}\,t_{j}^{-1}}{f_{j,1-\text{id}}\left(\mathbf{B}\right)}\\
\\
v_{10} & = & t_{j}\\
\\
v_{11} & = & s_{j}
\end{array},\end{cases}
\]
In order to satisfy Eq. (\ref{Slice dependence}) we may assume without
lost of generality that $t_{j}=1=s_{j}$, for all $0\le j<m$. Hence
we have 
\[
\begin{cases}
\begin{array}{ccc}
v_{00} & = & \frac{c_{j,1-\text{id}}}{f_{j,\text{id}}\left(\mathbf{B}\right)}\\
\\
v_{01} & = & \frac{c_{j,1-\text{id}}}{f_{j,1-\text{id}}\left(\mathbf{B}\right)}\\
\\
v_{10} & = & 1\\
\\
v_{11} & = & 1
\end{array}.\end{cases}
\]
Generically, the solution to system of equations
\[
\left\{ \prod_{0\le i<n}\mathbf{V}\left[i,\sigma\left(i\right)\right]=\sum_{\gamma\in\text{S}_{n}\backslash\left\{ \text{id}\right\} }\frac{c_{j,\gamma}\,\text{sgn}\left(\gamma\right)}{f_{j,\gamma}\left(\mathbf{B}\right)}\exp\left(\frac{2\pi\,\text{lex}\sigma\cdot\text{lex}\gamma}{n!}\sqrt{-1}\right)\right\} _{\sigma\in\text{S}_{n}}
\]
for distinct values of $j$ are distinct. Taking the solution for
some fixed index $j$ determine the entries of $\mathbf{V}$ and in
turn determine the entries of $\mathbf{M}$ in terms of the unknowns
$\left\{ c_{j,1-\text{id}}\right\} _{j}$. The requirement that 
\[
\mathbf{U}\left[:,1\right]=\mathbf{1}_{m\times1}.
\]
imposes additional constraints on the unknown $c_{j,1-\text{id}}$
which are generically inconsistent when $n=2$. If we take $m$ to
also be equal to $2$ then to resolve the consistency issue, it is
necessary and sufficient to require the solution $\mathbf{V}$ to
be the same for all index $j\in\left\{ 0,1\right\} $. This of course
is far from the generic setting and yields the constraints
\[
\det\left\{ \det\left(\mathbf{V}\right)\mathbf{I}_{2}-\mathbf{M}\left[:2,\,:2\right]\right\} =0=\det\left\{ \det\left(\mathbf{V}\right)\mathbf{I}_{2}-\mathbf{M}\left[2:4,\,2:4\right]\right\} 
\]
\[
\begin{cases}
\begin{array}{ccc}
\frac{\left(b_{001}b_{010}v_{01}v_{10}-b_{000}b_{011}v_{00}v_{11}\right)\left(v_{01}v_{10}-v_{00}v_{11}\right)\left(b_{000}-b_{001}\right)\left(b_{010}-b_{011}\right)}{b_{000}b_{001}b_{010}b_{011}} & = & 0\\
\\
\frac{\left(b_{101}b_{110}v_{01}v_{10}-b_{100}b_{111}v_{00}v_{11}\right)\left(v_{01}v_{10}-v_{00}v_{11}\right)\left(b_{100}-b_{101}\right)\left(b_{110}-b_{111}\right)}{b_{100}b_{101}b_{110}b_{111}} & = & 0
\end{array},\end{cases}
\]
Assuming that 
\[
\left(v_{01}v_{10}-v_{00}v_{11}\right)\ne0,\,\left(b_{000}-b_{001}\right)\left(b_{010}-b_{011}\right)\ne0,\,\left(b_{100}-b_{101}\right)\left(b_{110}-b_{111}\right)\ne0
\]
and 
\[
b_{000}b_{001}b_{010}b_{011}\ne0,\,b_{100}b_{101}b_{110}b_{111}\ne0
\]
then the system reduces to 
\[
\begin{cases}
\begin{array}{ccc}
b_{001}b_{010}v_{01}v_{10}-b_{000}b_{011}v_{00}v_{11} & = & 0\\
\\
b_{101}b_{110}v_{01}v_{10}-b_{100}b_{111}v_{00}v_{11} & = & 0
\end{array}, & \implies\left(b_{001}b_{010}b_{100}b_{111}-b_{101}b_{110}b_{000}b_{011}\right)=0\end{cases}
\]
which we recognize at once as hyperdeterminant of the $2\times2\times2$
as introduced in \cite{GNANG2017238}. However in the generic setting,
when $n=m>2$ the derivation results in more variables then there
are constraints since $n!>n$ whenever $n>2$, therefore generically
the constraints admit solutions from which the desired result follows.

\subsection{Necessary and sufficient conditions for the existence of inverse
pairs}

Inverse pairs, introduced in \cite{MB90,MB94}, extend to hypermatrices
the notion of matrix inverse and enable us to describe a hypermatrix
analog of the general linear group as follows 
\[
\text{GL}_{m\times n\times p}\left(m\times p\times p,\,p\times n\times n,\,\mathbb{K}\right)\,:=
\]
\begin{equation}
\left\{ \left(\mathbf{A},\mathbf{B}\right)\in\mathbb{K}^{m\times p\times p}\times\mathbb{K}^{p\times n\times p}:\,\exists\:\left(\mathbf{C},\mathbf{D}\right)\in\mathbb{K}^{m\times p\times p}\times\mathbb{K}^{p\times n\times p},\ \text{Prod}\left(\mathbf{C},\text{Prod}\left(\mathbf{A},\mathbf{X},\mathbf{B}\right),\mathbf{D}\right)=\mathbf{X},\ \forall\,\mathbf{X}\in\mathbb{K}^{m\times n\times p}\right\} \label{Inverse pair definition}
\end{equation}
By properties of the transpose it follows that 
\[
\left(\mathbf{A},\mathbf{B}\right)\in\text{GL}_{m\times n\times p}\left(m\times p\times p,\,p\times n\times n,\,\mathbb{K}\right),\;\exists\:\left(\mathbf{C},\mathbf{D}\right)\in\mathbb{K}^{m\times p\times p}\times\mathbb{K}^{p\times n\times p}
\]
such that
\[
\implies\begin{cases}
\begin{array}{ccc}
\text{Prod}\left(\text{Prod}\left(\mathbf{X}^{\top},\mathbf{B}^{\top},\mathbf{A}^{\top}\right),\mathbf{D}^{\top},\mathbf{C}^{\top}\right) & = & \mathbf{X}^{\top}\\
\\
\text{Prod}\left(\mathbf{D}^{\top^{2}},\mathbf{C}^{\top^{2}},\text{Prod}\left(\mathbf{B}^{\top^{2}},\mathbf{A}^{\top^{2}},\mathbf{X}^{\top^{2}}\right)\right) & = & \mathbf{X}^{\top^{2}}
\end{array} & \forall\:\mathbf{X}\in\mathbb{K}^{m\times n\times p}.\end{cases}
\]
Note that in Eq. (\ref{Inverse pair definition}) the hypermatrix
pair $\left(\mathbf{A},\mathbf{B}\right)$ is the inner-inverse pair
to of the pair $\left(\mathbf{C},\mathbf{D}\right)$. The simplest
illustration of a family of invertible pairs of hypermatrices is associated
with third order hypermatrix analog of the subgroup of scaling matrices
over a skew field $\mathbb{K}$. $\mathbf{A}\in\mathbb{K}^{m\times p\times p}$
and $\mathbf{B}\in\mathbb{K}^{p\times n\times p}$ having entries
specified by 
\[
\mathbf{A}\left[i,t,k\right]=\begin{cases}
\begin{array}{cc}
\alpha_{it}\in\mathbb{K}\backslash\left\{ 0\right\}  & \mbox{ if }0\le t=k<p\\
0 & \mbox{otherwise}
\end{array} & \forall\,\begin{cases}
\begin{array}{c}
0\le i<m\\
0\le t<p\\
0\le k<p
\end{array}\end{cases}\end{cases},
\]
\[
\mathbf{B}\left[t,j,k\right]=\begin{cases}
\begin{array}{cc}
\beta_{tj}\in\mathbb{K}\backslash\left\{ 0\right\}  & \mbox{ if }0\le t=k<p\\
0 & \mbox{otherwise}
\end{array} & \forall\,\begin{cases}
\begin{array}{c}
0\le t<p\\
0\le j<n\\
0\le k<p
\end{array}\end{cases}\end{cases},
\]
\[
\implies\mbox{Prod}\left(\mathbf{A},\mathbf{X},\mathbf{B}\right)\left[i,j,k\right]=\alpha_{ik}\,\mathbf{X}\left[i,j,k\right]\,\beta_{kj},\quad\forall\:\begin{cases}
\begin{array}{c}
0\le i<m\\
0\le j<n\\
0\le k<p
\end{array}\end{cases}.
\]
Consequently, the entries of the corresponding outer-inverse pair
$\left(\mathbf{C},\mathbf{D}\right)$ has entries specified by 
\[
\mathbf{C}\left[i,t,k\right]=\begin{cases}
\begin{array}{cc}
\alpha_{it}^{-1} & \mbox{ if }0\le t=k\le p\\
0 & \mbox{otherwise}
\end{array} & \forall\,\begin{cases}
\begin{array}{c}
0\le i<m\\
0\le t<p\\
0\le k<p
\end{array}\end{cases}\end{cases},
\]
\[
\mathbf{D}\left[t,j,k\right]=\begin{cases}
\begin{array}{cc}
\beta_{tj}^{-1} & \mbox{ if }0\le t=k\le p\\
0 & \mbox{otherwise}
\end{array} & \forall\,\begin{cases}
\begin{array}{c}
0\le t<p\\
0\le j<n\\
0\le k<p
\end{array}\end{cases}\end{cases}.
\]
We see that by analogy to the the matrix case the hypermatrix scaling
action forms a group over any skew field $\mathbb{K}$. Invertible
scaling hyprematrices therefore correspond to the largest subset of
invertible hypermatrix pairs whose inverse are obtained by inverting
the non-zero entries in the given pair. Note that in general the third
order analog of the general linear group does not form a group. Moreover
entry-wise Eq. (\ref{Inverse pair definition}) yields the constraints
\begin{equation}
\mathbf{X}\left[i,j,k\right]=\sum_{0\le t<p}\mathbf{C}\left[i,t,k\right]\,\mbox{Prod}\left(\mathbf{A},\mathbf{X},\mathbf{B}\right)\left[i,j,t\right]\,\mathbf{D}\left[t,j,k\right],\:\forall\:\begin{cases}
\begin{array}{c}
0\le i<m\\
0\le j<n\\
0\le k<p
\end{array}\end{cases}.\label{condition 1}
\end{equation}
By substituting 
\[
\mbox{Prod}\left(\mathbf{A},\,\mathbf{X},\,\mathbf{B}\right)\left[i,j,t\right]=\sum_{0\le s<p}\mathbf{A}\left[i,s,t\right]\,\mathbf{X}\left[i,j,s\right]\,\mathbf{B}\left[s,j,t\right]
\]
into (\ref{condition 1}) we get 
\[
\mathbf{X}\left[i,j,k\right]=\sum_{0\le s,t<p}\mathbf{C}\left[i,t,k\right]\,\mathbf{A}\left[i,s,t\right]\,\mathbf{X}\left[i,j,s\right]\,\mathbf{B}\left[s,j,t\right]\,\mathbf{D}\left[t,j,k\right],\:\forall\begin{cases}
\begin{array}{c}
0\le i<m\\
0\le j<n\\
0\le k<p
\end{array}\end{cases}.
\]
Over skew fields it is not at all clear how to solve in general such
systems. However if we are instead working over commutative algebraically
closed field such as $\mathbb{C}$, We may regroup and reorder the
factors in the summands as follows
\[
\mathbf{X}\left[i,j,k\right]=\sum_{0\le s,t<p}\left(\mathbf{C}\left[i,t,k\right]\,\mathbf{D}\left[t,j,k\right]\right)\,\left(\mathbf{A}\left[i,s,t\right]\,\mathbf{B}\left[s,j,t\right]\right)\,\mathbf{X}\left[i,j,s\right],\:\forall\:\begin{cases}
\begin{array}{c}
0\le i<m\\
0\le j<n\\
0\le k<p
\end{array}\end{cases}.
\]
Consequently, the inverse pair relation (\ref{Inverse pair definition})
asserts that 
\[
\forall\:\begin{array}{c}
0\le i<m\\
0\le j<n\\
0\le k<p
\end{array},\:\sum_{0\le t<p}\left(\mathbf{C}\left[i,t,k\right]\,\mathbf{D}\left[t,j,k\right]\right)\left(\mathbf{A}\left[i,s,t\right]\,\mathbf{B}\left[s,j,t\right]\right)=\begin{cases}
\begin{array}{cc}
1 & \mbox{ if }0\le k=s<p\\
0 & \mbox{otherwise}
\end{array} & .\end{cases}
\]
This in effect flattens the inverse pair relation to finding elements
of the matrix General linear group GL$_{m\cdot n\cdot p}\left(\mathbb{C}\right)$.
Let $\mathbf{F}=\mathcal{F}\left(\mathbf{A},\mathbf{B}\right)\in$
GL$_{m\cdot n\cdot p}\left(\mathbb{C}\right)$ expressed in terms
of the entries of $\mathbf{A}$ and $\mathbf{B}$ as follows:
\[
\forall\,0\le i<m,\,0\le j<n,\,0\le s,t<p,
\]
\[
\mathbf{F}\left[i\cdot n\cdot p+j\cdot p+t,\,i\cdot n\cdot p+j\cdot p+s\right]=\mathbf{A}\left[i,s,t\right]\cdot\mathbf{B}\left[s,j,t\right].
\]
Note that $\mathbf{F}$ is a direct sum of $m\cdot n$ elements of
GL$_{p^{2}}\left(\mathbb{C}\right)$ . Assuming the pair $\mathbf{A}\in\mathbb{C}^{m\times p\times p}$,
$\mathbf{B}\in\mathbb{C}^{p\times n\times p}$ given, we determine
the corresponding inverse pair $\mathbf{C}\in\mathbb{C}^{m\times p\times p}$,
$\mathbf{D}\in\mathbb{C}^{p\times n\times p}$. Given that we know
the entries of $\mathbf{F}$ in terms of the entries of $\mathbf{A}$
and $\mathbf{B}$, it follows that 
\[
\mathbf{F}^{-1}=\frac{\mbox{Adj}\left(\mathbf{F}\right)}{\det\left(\mathbf{F}\right)}.
\]
Consequently the entries of $\mathbf{C}$ and $\mathbf{D}$ are determined
by solving system type 2 as described in \cite{2018arXiv180803743G}
given by 
\begin{equation}
\mathcal{S}:=\left\{ \mathbf{C}\left[i,t,k\right]\cdot\mathbf{D}\left[t,j,k\right]=\mathbf{F}^{-1}\left[i\cdot n\cdot p+j\cdot p+k,\,i\cdot n\cdot p+j\cdot p+t\right]\right\} _{\begin{array}{c}
0\le i<m,\,0\le j<n\\
0\le k,t<p
\end{array}}.\label{Condition 2}
\end{equation}
Therefore the two necessary and sufficient conditions which guarantee
that the hypermatrix pair $\mathbf{A}\in\mathbb{C}^{m\times p\times p}$,
$\mathbf{B}\in\mathbb{C}^{p\times n\times p}$ admit an outer inverse
pair $\mathbf{C}\in\mathbb{C}^{m\times p\times p}$, $\mathbf{D}\in\mathbb{C}^{p\times n\times p}$
are that det$\left(\mathbf{F}\right)\neq0$ and that $\mathbf{F}^{-1}$
lies in the log linear column space of the exponent matrix of the
system $\mathcal{S}$ as discussed in \cite{2018arXiv180803743G}.

\subsection{Third order hypermatrix rank\textendash nullity theorem}

The third order hypermatrix variant of the rank\textendash nullity
theorem generalizes the well-known rank\textendash nullity theorem
of linear algebra \cite{10.2307/2324660}. For the sake of completeness,
we recall the statement and a proof of the rank\textendash nullity
theorem. We subsequently extend the argument to third order hypermatrices.\\
\\
\textbf{Theorem 8} : $\mathbf{A}\in\mathbb{C}^{m\times n}$ ( where
$r\le n=$ min$\left(m,n\right)$ ) has nullity $\left(n-r\right)$
iff Rank$\left(\mathbf{A}\right)=r$.\\
\\
\emph{Proof of sufficiency} : It follows from the hypothesis that
there exist an invertible matrix $\mathbf{X}$ of size $n\times n$,
such that the last $\left(n-r\right)$ columns of $\mathbf{A}\mathbf{X}$
are zero columns. Hence 
\[
\text{Prod}\left(\text{Prod}\left(\mathbf{A},\mathbf{X}\right),\mathbf{X}^{-1}\right)=\mathbf{A}=\sum_{0\le t<r}\mbox{Prod}_{\boldsymbol{\Delta}^{(t)}}\left(\text{Prod}\left(\mathbf{A},\mathbf{X}\right),\mathbf{X}^{-1}\right)
\]
thereby expressing $\mathbf{A}$ as a sum of $r$ outer products.\emph{}\\
\emph{}\\
\emph{Proof of necessity} : Assuming that Rank$\left(\mathbf{A}\right)=r$,
we exhibit $n-r$ columns of an invertible matrix which are mapped
to zero by $\mathbf{A}$ 
\begin{equation}
\mbox{Rank}\left(\mathbf{A}\right)=r\Rightarrow\begin{cases}
\begin{array}{c}
\exists\,S\subset\left\{ 0,\cdots,n-1\right\} \\
\left|S\right|=r
\end{array}\end{cases}\:\text{ such that }\:\mathbf{A}=\sum_{t\in S}\text{Prod}_{\boldsymbol{\Delta}^{(t)}}\left(\mathbf{U},\mathbf{V}\right),\label{Rank criteria}
\end{equation}
For some $\mathbf{U}\in\mathbb{C}^{m\times n}$ and $\mathbf{V}\in\mathbb{C}^{n\times n}$.
Let us replace by zero columns the columns of $\mathbf{U}$ indexed
by $\left\{ 0,\cdots,m-1\right\} -S$ and also replace by zero rows
the rows of $\mathbf{V}$ indexed by $\left\{ 0,\cdots,n-1\right\} -S$.
This substitution leaves the identity (\ref{Rank criteria}) unchanged.
Note that the rows of $\mathbf{V}$ must be linearly independent for
if this was not case by Theorem 6 the Rank$\left(\mathbf{A}\right)$
would be less then $r$. Consequently, we may replace the rows of
$\mathbf{V}$ indexed by $\left\{ 0,\cdots,n-1\right\} -S$ by an
arbitrary basis of the orthogonal complement of the space spanned
by the rows of $\mathbf{V}$ indexed by $S$. Such a basis can always
be taken from some particular subset of rows of the identity matrix.
So long as the columns of $\mathbf{U}$ indexed by $\left\{ 0,\cdots,m-1\right\} -S$
remain zero columns the row substitutions in $\mathbf{V}$ leave the
identity (\ref{Rank criteria}) unchanged. Consequently det$\left(\mathbf{V}\right)\ne0$
and the columns of $\mathbf{A}\cdot\mathbf{V}^{-1}$ indexed by $\left\{ 0,\cdots,n-1\right\} -S$
are mapped by $\mathbf{A}$ to the zero columns thus concluding our
proof.\\
\\
\\
The nullity of a hypermatrix $\mathbf{A}\in\mathbb{C}^{m\times n\times p}$
is defined to be the maximum number distinct pairs of depth slices
of an invertible pair which are mapped zero slices by the action of
an invertible pair of hypermatrices on an input hypermatrix $\mathbf{A}$.\\
\\
\\
\textbf{Theorem 9} : $\mathbf{A}\in\mathbb{C}^{m\times n\times p}$
(where $r\le p=$min$\left\{ m,n,p\right\} $) has nullity $\left(p-r\right)$
iff Rank$\left(\mathbf{A}\right)=r$.\\
\\
\emph{Proof sufficiency} : It follows from the hypothesis that there
exists an invertible pair of hypermatrices $\mathbf{X}_{0}$, $\mathbf{X}_{1}$,
respectively of size $m\times p\times p$ and $p\times n\times p$,
such that the last $\left(p-r\right)$ depth slices of $\mbox{Prod}\left(\mathbf{X}_{0},\mathbf{A},\mathbf{X}_{1}\right)$
are zero depth slices.\\
Let $\left(\mathbf{Y}_{0},\mathbf{Y}_{1}\right)$, respectively of
size $m\times p\times p$ and $p\times n\times p$ denote the outer-inverse
pair of $\left(\mathbf{X}_{0},\mathbf{X}_{1}\right)$, hence 
\[
\mbox{Prod}\left(\mathbf{Y}_{0},\mbox{Prod}\left(\mathbf{X}_{0},\mathbf{A},\mathbf{X}_{1}\right),\mathbf{Y}_{1}\right)=\mathbf{A}=\sum_{0\le t<r}\mbox{Prod}_{\boldsymbol{\Delta}^{(t)}}\left(\mathbf{Y}_{0},\mbox{Prod}\left(\mathbf{X}_{0},\mathbf{A},\mathbf{X}_{1}\right),\mathbf{Y}_{1}\right)
\]
thereby expressing $\mathbf{A}$ as a sum of $r$ outer products.\\
\\
\emph{Proof of necessity} : Assuming that Rank$\left(\mathbf{A}\right)=r$,
we exhibit $\left(p-r\right)$ pairs of depth slices of an invertible
hypermatrix pair which are mapped to zero by $\mathbf{A}$ 
\begin{equation}
\mbox{Rank}\left(\mathbf{A}\right)=r\Rightarrow\begin{array}{c}
\exists\,S\subset\left\{ 0,\cdots,n-1\right\} \\
\left|S\right|=r
\end{array}\:\text{ such that }\:\mathbf{A}=\sum_{t\in S}\text{Prod}_{\boldsymbol{\Delta}^{(t)}}\left(\mathbf{U},\mathbf{V},\mathbf{W}\right),\label{Hypermatrix rank criteria}
\end{equation}
For some $\mathbf{U}\in\mathbb{C}^{m\times p\times p}$ , $\mathbf{V}\in\mathbb{C}^{m\times n\times p}$
and $\mathbf{W}\in\mathbb{C}^{p\times n\times p}$. Let us replace
by zero column slices the column slices of $\mathbf{U}$ indexed by
$\left\{ 0,\cdots,m-1\right\} -S$, also replace by zero depth slices
of $\mathbf{V}$ the depth slices of $\mathbf{V}$ indexed by $\left\{ 0,\cdots,n-1\right\} -S$
and finally replacing by zero row slices of $\mathbf{W}$ the row
slices of $\mathbf{W}$ indexed by $\left\{ 0,\cdots,n-1\right\} -S$.
This substitution leaves the identity ( \ref{Hypermatrix rank criteria}
) unchanged. Note that the depth slices of $\mathbf{V}$ do not satisfy
the assumptions of Theorem 7 otherwise by Theorem 7 the Rank$\left(\mathbf{A}\right)$
would be less then $r$. Consequently, we may replace the column slices
of $\mathbf{U}$ indexed by $\left\{ 0,\cdots,n-1\right\} -S$ and
also replace the row slices of $\mathbf{W}$ indexed by $\left\{ 0,\cdots,n-1\right\} -S$
so as to ensure that $\mathbf{U}$, $\mathbf{V}$ form an invertible
pair. So long as the depth slices of $\mathbf{V}$ indexed by $\left\{ 0,\cdots,m-1\right\} -S$
remain zero slices, the substitutions into $\mathbf{U}$ and $\mathbf{V}$
leave the identity ( \ref{Hypermatrix rank criteria} ) unchanged.
The depth slices of Prod$\left(\mathbf{U}^{-1},\mathbf{A},\mathbf{W}^{-1}\right)$
indexed by $\left\{ 0,\cdots,n-1\right\} -S$ are mapped by $\mathbf{A}$
to the zero depth slices, thus concluding our proof.

\bibliographystyle{amsalpha}
\bibliography{mybib}

\newcommand{\etalchar}[1]{$^{#1}$}
\providecommand{\bysame}{\leavevmode\hbox to3em{\hrulefill}\thinspace}
\providecommand{\MR}{\relax\ifhmode\unskip\space\fi MR }
\providecommand{\MRhref}[2]{%
  \href{http://www.ams.org/mathscinet-getitem?mr=#1}{#2}
}
\providecommand{\href}[2]{#2}
\begin{thebibliography}{BCC{\etalchar{+}}16}

\bibitem[BCC{\etalchar{+}}16]{2016arXiv160506702B}
J.~{Blasiak}, T.~{Church}, H.~{Cohn}, J.~A. {Grochow}, E.~{Naslund}, W.~F.
  {Sawin}, and C.~{Umans}, \emph{{On cap sets and the group-theoretic approach
  to matrix multiplication}}, ArXiv e-prints (2016).

\bibitem[GER11]{GER}
E.~K. Gnang, A.~Elgammal, and V.~Retakh, \emph{A spectral theory for tensors},
  Annales de la faculte des sciences de Toulouse Mathematiques \textbf{20}
  (2011), no.~4, 801--841.

\bibitem[GF17a]{2017arXiv170606090G}
E.~K. {Gnang} and Y.~{Filmus}, \emph{{On the Bhattacharya-Mesner rank of third
  order hypermatrices}}, ArXiv e-prints (2017).

\bibitem[GF17b]{Gnang2017238}
Edinah~K. Gnang and Yuval Filmus, \emph{On the spectra of hypermatrix direct
  sum and kronecker products constructions}, Linear Algebra and its
  Applications \textbf{519} (2017), 238 -- 277.

\bibitem[GG18]{2018arXiv180803743G}
E.~K. {Gnang} and J.~S. {Gnang}, \emph{{Sketch for a Theory of Constructs}},
  ArXiv e-prints (2018).

\bibitem[GKZ94]{GKZ}
I.~Gelfand, M.~Kapranov, and A.~Zelevinsky, \emph{Discriminants, resultants and
  multidimensional determinant}, Birkhauser, Boston, 1994.

\bibitem[GR91]{GR}
I.~M. Gelfand and V.~S. Retakh, \emph{Determinants of matrices over
  noncommutative rings}, Functional Analysis and Its Applications \textbf{25}
  (1991), no.~2, 91--102.

\bibitem[GR97]{Gelfand1997}
I.~Gelfand and V.~Retakh, \emph{Quasideterminants, i}, Selecta Mathematica
  \textbf{3} (1997), no.~4, 517--546.

\bibitem[GRW05]{GRW}
I.~M. Gelfand, V.~S. Retakh, and R.~Wilson, \emph{Quasideterminants}, Advances
  in Mathematics \textbf{193} (2005), no.~1, 56--141.

\bibitem[Ker08]{RK}
Richard Kerner, \emph{Ternary and non-associative structures}, International
  Journal of Geometric Methods in Modern Physics \textbf{5} (2008), 1265--1294.

\bibitem[Lim13]{Lim2013}
Lek-Heng Lim, \emph{Tensors and hypermatrices}, Handbook of Linear Algebra
  (Leslie Hogben, ed.), CRC Press, 2013.

\bibitem[MB90]{MB90}
D.~M. Mesner and P.~Bhattacharya, \emph{Association schemes on triples and a
  ternary algebra}, Journal of combinatorial theory \textbf{A55} (1990),
  204--234.

\bibitem[MB94]{MB94}
D.~M. Mesner and P.~Bhattacharya, \emph{A ternary algebra arising from
  association schemes on triples}, Journal of algebra \textbf{164} (1994),
  595--613.

\bibitem[S{\etalchar{+}}15]{sage}
W.\thinspace{}A. Stein et~al., \emph{{S}age {M}athematics {S}oftware ({V}ersion
  6.9)}, The Sage Development Team, 2015, {\tt http://www.sagemath.org}.

\bibitem[Str93]{10.2307/2324660}
Gilbert Strang, \emph{The fundamental theorem of linear algebra}, The American
  Mathematical Monthly \textbf{100} (1993), no.~9, 848--855.

\end{thebibliography}

\end{document}